\documentclass[10pt]{article}
\usepackage{amsfonts}
\def\bea{\begin{eqnarray}}
\def\eea{\end{eqnarray}}
\def\beann{\begin{eqnarray*}}
\def\eeann{\end{eqnarray*}}
\def\beq{\begin{equation}}
\def\eeq{\end{equation}}
\def\ba{\begin{array}}
\def\ea{\end{array}}
\def\ben{\begin{enumerate}}
\def\een{\end{enumerate}}
\def\f{\phi}
\def\m{\mu}
\def\pd{\partial_}
\def\p{\rho}

\def\la{\lambda}
\def\a{\alpha}

\def\s{\sigma}
\def\om{\omega}

\newtheorem{theorem}{Theorem}
\newtheorem{lemma}{Lemma}
\newtheorem{proposition}{Proposition}
\newtheorem{corollary}{Corollary}
\newtheorem{remark}{Remark}
\newtheorem{con}{Conjecture}

\font\mybb=msbm10 at 11pt

\def\bb#1{\hbox{\mybb#1}}

\def\bZ {\bb{Z}}
\def\bR {\bb{R}}

\def\bH {\bb{H}}
\def\bC {\bb{C}}

\def\bS {\bb{S}}
\def\bO {\bb{O}}
\def\e  {\epsilon}

\def\dd{d^\dagger}
\def\cG{{\buildrel {\it c} \over \Gamma}}

\def\cL{{\cal L}}

\date{}
\begin{document}
\thispagestyle{empty}
\begin{center}
\vspace {.7cm} {\large {\bf{Deformations of
 generalized calibrations and compact
 non-K{\"a}hler manifolds with
vanishing first Chern class}} }
 \vskip 1.0truecm
 \vskip 1.0truecm
  {\large{\bf{
Jan Gutowski${}^1$, Stefan Ivanov${}^2$ and George
  Papadopoulos${}^3$}}}
 \vskip 2.0truecm
 {\normalsize{\sl 1. Department of Physics}}
\\ {\normalsize{\sl Queen Mary University of London}}
\\ {\normalsize{\sl Mile End, London E1 4NS.}}
\vskip 1.0truecm
 {\normalsize{\sl 2. Department of Mathematics}}
\\ {\normalsize{\sl University of Sofia}}
\\ {\normalsize{\sl \lq\lq St. Kl.. Ohridski''}}
\vskip 1.0truecm {\normalsize{\sl 3. Department of Mathematics}}
\\ {\normalsize{\sl King's College London}}
\\ {\normalsize{\sl Strand}}
\\ {\normalsize{\sl London WC2R 2LS}}
\vskip 2.0truecm
{\bf {Abstract}}
\end{center}
We  investigate the deformation theory of a class of generalized
calibrations in Riemannian manifolds for which the tangent bundle
has reduced structure group $U(n)$, $SU(n)$, $G_2$ and $Spin(7)$.
{}For this we use the property of the associated calibration form
to be parallel with respect to a metric connection which may have
non-vanishing torsion. In all these cases, we find that if there
is a moduli space, then it is finite dimensional.
 We present various examples
of generalized calibrations that include almost hermitian
manifolds with structure group $U(n)$ or $SU(n)$, nearly parallel
$G_2$ manifolds and group manifolds. We find that some Hopf
fibrations are deformation families of generalized calibrations.
In addition, we give sufficient conditions for a hermitian
manifold $(M,g,J)$ to admit Chern and Bismut connections with
holonomy contained in $SU(n)$. In particular we show that any
connected sum of $k \geq 3$ copies of $S^3 \times S^3$
admits a hermitian structure for which the restricted
holonomy of a Bismut connection is contained in $SU(3)$.

\section{Introduction}

 Riemannian manifolds with structure
group  a subgroup of an orthogonal group under mild topological
assumptions admit a  connection for which its reduced holonomy is a
subgroup of the structure group. This connection   is not
necessarily the Levi-Civita connection but it may have
non-vanishing torsion. The existence of such a connection with
reduced holonomy a subgroup of an orthogonal group
 does not imply other  geometric properties
on a Riemannian manifold, like for example irreducibility.
This is unlike the well-known
  case that involves the reduction of the holonomy group of
the Levi-Civita connection which has led to  Berger's
classification list. Nevertheless the question arises as to
whether the reduction of the structure group of a Riemannian manifold is
related to some underlying geometric structure.

The aim of this paper is three-fold. First, we shall show
 that Riemannian manifolds  which admit a metric connection
 with holonomy an appropriate subgroup
of the orthogonal  group may have submanifolds which are
calibrated with respect to a generalized calibration. Second we
shall investigate the moduli space of these calibrated
submanifolds. Finally, we shall show the existence of a large
class of hermitian
 manifolds with trivial canonical bundle
which admit either a Chern or a Bismut connection which has
 reduced holonomy contained in $SU(n)$.
Our latter result can be though off as a generalization of the
Calabi-Yau theorem in the context
of hermitian manifolds which are not K{\"a}hler.

Generalized calibrations were introduced by  Gutowski and
 Papadopoulos \cite{gutpap} and further investigated in
\cite{gutpapt} to describe the solitons of
brane actions with a non-vanishing Wess-Zumino term. These
solitons are certain  submanifolds  which minimize an energy
functional and are associated with calibration forms.  These
forms, unlike the case of standard calibrations,
 are {\it not} closed.  In what follows we shall use the
term {\it generalized calibration} to refer to both the
calibration form and the calibrated submanifold. The distinction
between the two will be clear from the context.

In this paper, we shall  demonstrate that generalized calibrations
arise in the investigation of manifolds that
 admit a metric connection with possibly
{\it non-vanishing torsion} which has holonomy an appropriate
subgroup of the orthogonal group. Although, generalized
calibrations can be investigated independently, their use in
manifold theory becomes more transparent in the context of
holonomy groups and in
 manifolds with reduced structure group.
 This is because for certain holonomy groups, like for example those
that occur in Berger's list; $U(n)$ (2n), $SU(n)$ (2n),
$Sp(n)\cdot Sp(1)$ (4n),
 $Sp(n)$(4n), $G_2$ (7) and $Spin(7)$ (8),
manifolds admit {\it parallel} calibration forms which however are
{\it not} necessarily closed; in parenthesis we have denoted the
real dimension of the associated manifolds. Such forms give rise
to generalized calibrated submanifolds which are minima of the
energy functional
\beq
 E(Z)={\rm Vol}(Z)-\int_Z \psi\ ,
\label{stuffa}
\eeq
where $Z$
is a k-dimensional submanifold and $\psi$ is  a calibration form
of degree $k$, $d\psi\not=0$.  The submanifolds that minimize $E$
 are {\it not} necessarily   minimal.

We shall focus our investigation to the generalized calibrations
associated with the holonomy groups $U(n)$ (2n), $SU(n)$ (2n),
$G_2$ (7) and $Spin(7)$ (8). We shall show that in most of these
cases, the differential system associated with the deformation of
the above generalized calibrations is elliptic. So if the moduli
space of a generalized calibration exists, then it is finite
dimensional. We shall not investigate the obstruction theory; this
will appear in another publication. We shall also compute the
second variation of the energy functional. The differential
systems that arise in the deformation of generalized calibrations
will also be investigated for various classes of manifolds that
admit connections with the above holonomy. We shall see that in
some cases they become simplified. We shall also
 give a large number of
generalized calibrations as submanifolds of group manifolds,
complex manifolds and  homogeneous spaces. In particular, we shall
show that some {\it Hopf fibrations} are {\it families} of
generalized calibrations.

In the second part of the paper we shall  focus on hermitian
manifolds with vanishing first Chern Class. This is because they
are a generalization of Calabi-Yau manifolds. Recently such
manifolds have found applications in the investigation of Reid's
conjecture and of mirror symmetry. This conjecture can be stated
as follows: Let $X$ be a three-dimensional
 Calabi-Yau manifold and suppose that $X$ can
 be blown down along a rational
 curve to a possibly singular manifold $Y_1$. Such
 singularities of Calabi-Yau
 manifolds can be removed by a small deformation.
 Let $\tilde Y_1$ be the smooth
 deformation of $Y$. Now $Y_1$ has trivial canonical
 bundle and $b_2(\tilde Y_1)=b_2(X)-1$.
 Continuing this procedure, we shall end up with a
 smooth manifold $\tilde X$
 with trivial canonical bundle and $b_2(\tilde X)=0$.
 So $\tilde X$ cannot be K{\"a}hler.
 The conjecture is that if $X$ and $Z$ are Calabi-Yau
 manifolds with $b_3(X)=b_3(Z)$,
 then $\tilde Z$ is in the same deformation class of $\tilde X$.

The canonical bundle of a hermitian manifold can be {\it
topologically} but not {\it holomorphically} trivial. For the
definition of the former we take that the first Chern class
vanishes. For the definition of the latter, we take that the
canonical bundle admits a nowhere vanishing holomorphic section.
We remark that there are canonical bundles which are topologically
but not holomorphically trivial, such as the canonical bundle of
$SU(3)$. It was discovered recently by Hitchin \cite{Hit} that
complex three-folds with holomorphically trivial canonical bundle
appear as critical points of a certain diffeomorphism invariant
functional on the space of differential three-forms on a closed
six-dimensional manifold.

There are several connections on hermitian manifolds compatible
with both the hermitian metric and complex structure which
coincide with the Levi-Civita connection in the K{\"a}hler case.
Amongst these connections the Chern connection is the unique
connection for which the torsion 2-form is of type (2,0)+(0,2),
 and the
Bismut connection for which the torsion is a three form. The
latter connection was used by Bismut \cite{Bis} to prove a local
index formula for the Dolbeault operator when the manifold is not
K{\"a}hler; vanishing theorems for the Dolbeault cohomology on a
compact Hermitian non-K{\"a}hler manifold were found
\cite{AI1,IP,IP1}. For other applications of the Bismut connection
see \cite{Stro,IP,IP1}. In particular in \cite{IP1} obstructions
have been found to the Hodge numbers $h^{0,1}, h^{0,n}$ for
hermitian manifolds whose Bismut connection has reduced holonomy
contained in $SU(3)$.

Given a hermitian manifold $(M,g,J)$ with vanishing first Chern
class, we give some sufficient conditions for $(M,g,J)$ to admit a
Chern or a Bismut connection with restricted holonomy contained in
$SU(n)$. The main tool that we shall use for the investigation of
hermitian manifolds with trivial canonical bundle is the
$\partial\bar\partial$-lemma. This lemma is valid for any compact
K{\"a}hler
 manifold but there are non-K{\"a}hler spaces
satisfying the $\partial\bar\partial$-lemma. A result of Deligne
states that any Moishezon  manifold is cohomologically K{\"a}hler
and therefore it satisfies the $\partial\bar\partial$-lemma. The
$\partial\bar\partial$-lemma also holds for any compact
non-K{\"a}hler 3-fold with holomorphically trivial canonical
bundle which is diffeomorphic to connected sums of $k\geq
2$-copies of $S^3\times S^3$ \cite{LuT}. One of our main goals is
to prove the following
\begin{theorem}
\label{t1} On a connected sum of $k\geq 2$-copies of $S^3\times
S^3$ there exists a hermitian structure with restricted holonomy
of the Bismut connection contained in $SU(3)$.
\end{theorem}

The celebrated Yau's solution of the Calabi conjecture \cite{Yau}
states that on a 2n-dimensional compact complex manifold with
vanishing first Chern class of K{\"a}hler type there exists a
K{\"a}hler metric with restricted holonomy contained in SU(n)
(Ricci flat K{\"a}hler metric). For  non-K{\"a}hler manifolds, it
appears that the following holds:
\begin{con}
 On any 2n-dimensional compact complex manifold with
vanishing first Chern class there exists a hermitian structure
with restricted holonomy of the Bismut connection contained in
SU(n).
\end{con}

Clearly this statement is true for connected sums of $k\geq
2$-copies of $S^3\times S^3$ in view of Theorem~\ref{t1}. It is
also true for  Moishezon manifolds and for compact complex
manifolds with vanishing first Chern class which are
cohomologically K\"ahler as we demonstrate in sections 17 and 18
below.

This paper has been organized as follows: In sections two and
three, we give the definition of generalized calibrations and
introduce the energy functional. In section four, we compute the
second variation of the energy functional and demonstrate the
relation between generalized calibrations and reduced holonomy. In
section five, we examine the deformations of a class of almost
hermitian calibrations. In section six, we derive the deformation
equations of SAS calibrations and in section seven, we explore
them in various special cases. In section eight, we give many
examples of SAS calibrations and deformation families. In section
nine, we give the deformation equations of generalized
co-associative calibrations. In section ten, we give the
deformation equations of generalized associative calibrations. In
section eleven, we investigate the deformation equations of
generalized associative and co-associative calibrations in various
manifolds with special $G_2$ structures.  In section twelve, we
give many examples of associative and co-associative calibrations
that include various deformation families. Some Hopf fibrations
are such deformation families. In section thirteen, we give the
deformation equations of generalized Cayley calibrations and in
section fourteen we give a group manifold example. In section
fifteen, we summarize some useful formulae for hermitian
manifolds. In section sixteen, we investigate the existence of
Chern connections with holonomy $SU(n)$ on hermitian manifolds
with trivial canonical bundles and in section seventeen we
investigate the existence of Bismut connections with holonomy
$SU(n)$ on hermitian manifolds with trivial canonical bundles.
In section eighteen, we give the proof of theorem one and in section
nineteen, we give examples of manifolds with the holonomy of Bismut
connection contained in $SU(3)$.

{\bf Acknowledgements.} We would like to thank T.Pantev and V.
 Tsanov for
helpful discussions.
 J.G. is supported by an
EPSRC postdoctoral grant.
S.I. is partially supported by Contract MM
809/1998 with the Ministry of Science and Education of Bulgaria,
Contract 353/2000 with the University of Sofia ``St. Kl. Ohridski".
S.I. is a member of the EDGE, Research Training Network
HPRN-CT-2000-00101, supported by the European  Human Potential
Programme.
G.P. is supported by Royal Society University
Research fellowship. Part of this work was done while one of us
G.P. was participating at the M-theory programme of the Newton
Institute.

\section{Generalized calibrations and relative de Rham Cohomology}

Here we shall describe some of the main properties of generalized
 calibrations using relative de Rham cohomology.
Generalized calibrations were  defined in
 \cite{gutpap} in the context of understanding the solitons
of brane actions with a Wess-Zumino term. In manifold theory,
 these solitons
are  certain submanifolds of a manifold which admits an
 appropriate form.
These submanifolds are not minimal but they are the minima of
a certain energy functional. We begin with a definition of
 generalized calibrations
and in particular of the energy functional.

\vskip 0.3cm
\leftline{{\bf Definition}}

 A generalized
calibration of degree $k$ is a k-form $\phi$ on an oriented manifold
$M$ which satisfies at every point $p$ the inequality
$\phi(\xi)|_p\leq 1$ for every oriented k-plane $\xi$ in $T_pM$.
\vskip 0.3cm

For {\it standard} calibrations  it is assumed in addition
that $\phi$ is
{\it closed}, $d\phi=0$. This is not the case here.

\vskip 0.3cm
\leftline{{\bf Definition}}
The contact set $C_p(\phi)$ at a point $p\in M$ of a
 calibration $\phi$ is
\beq
C_p=\{\xi\in Gr(k, T_pM):~ \phi(\xi)=1\} \ .
\label{stuffb}
\eeq
\vskip 0.3cm

For calibrations of interest the contact sets are not empty.
\vskip 0.3cm \leftline{{\bf Definition}} Generalized calibrated
submanifolds $X$ of $M$ are those for which $\phi(T_pX)=1$ at
every $p\in X$. \vskip 0.3cm In what follows we shall refer to
both $\phi$ and $X$ as ``generalized calibrations'' of $M$. The
distinction between the two will be clear from the context.

Generalized calibrations
 minimize a family of functionals
\cite{gutpap}. Here we shall repeat the analysis  using relative
homology and relative de Rham cohomology. Let $M$ be a manifold
and $N$ a submanifold of $M$ with ${\rm dim}N\geq k$. Suppose that
$(\alpha, \beta)$ is a pair of forms such that $\alpha\in
\Omega^{k+1}(M)$ and $\beta\in\Omega^k(N)$.
 In addition assume that $\beta$ is a
calibration. Next suppose that $K$ is a submanifold of $N$. Take an
open ball $D\subset K$ and consider the
functional
\begin{equation}\label{mfun}
{\cal E}(D, L)={\rm Vol}(D) -\int_{L}\alpha\ ,
\end{equation}
where $L$ is a submanifold of $M$ such that $\partial L=D+Y$. We
shall refer to ${\cal E}$ as the energy of $D$.

We shall show that the  functional ${\cal E}$ is minimized whenever
$D\subset X$  and $X$ is a calibrated submanifold of
$N$. However before we proceed to show this consider $D_1$ and
$D_2$ two open balls in $N$ such that $\partial
D_1=\partial D_2$. Then we have
\begin{equation}
\int_{D_1}\beta -\int_{L_1}\alpha-\big(\int_{D_2}\beta-
\int_{L_2}\alpha\big)
= \int_{S}\beta-\int_Z\alpha
\end{equation}
where $\partial L_2 =Y+D_2$ and $S=D_1-D_2$ is the sphere in $N$
 which can be constructed by gluing the discs
$D_1$ and $D_2$ along the common boundary taking into account their
 relative orientations. In addition $Z$ is
obtained by gluing $L_1$ and $L_2$ along $Y$ and so it has
 boundary $\partial Z=\partial L_1-\partial L_2=
(Y+D_1)-(Y+D_2)=D_1-D_2=S$. So $Z$ is a cycle in $M$ relative to
 the submanifold $N$ and $Z\in H_{k+1}(M,N)$.

Now suppose that $(\alpha, \beta)$ represents a trivial class
in the relative de Rham cohomology
$H^{k+1}_{dR}(M,N)$. Recall that the cohomology operator $d$
in relative de Rham cohomology is defined as
$d(\alpha, \beta)=(d\alpha, \alpha|_N-d\beta)$.
Therefore $[(\alpha, \beta)]$ is
a trivial class iff
$\alpha$ is  exact, $\alpha=d\gamma$, which implies that
 $d(\gamma|_N-\beta)=0$, and  $\gamma|_N-\beta$
 is an exact form in $N$, $\gamma|_N-\beta=d\zeta$.
In such a case, we have
\begin{equation}
\int_{S}\beta-\int_Z\alpha=0
\end{equation}
and so
\begin{equation}
\int_{D_1}\beta -\int_{L_1}\alpha-\big(\int_{D_2}\beta
-\int_{L_2}\alpha\big)=0\ .
\end{equation}
\begin{theorem}
Let $(\alpha, \beta)$ represent the trivial class in
$H^{k+1}_{dR}(M,N)$ and $\beta$ be a calibration form in $N$. Then
calibrated submanifolds $X$ of $N$ minimize the functional ${\cal
E}(D_1, L_1)$.
\end{theorem}
{\it Proof:} Let $X$ be a calibrated submanifold of $N$ and $D_1$
a disc in $X$. Then we have
\beq
 {\cal E}(D_1, L_1)={\rm
Vol}(D_1)-\int_{L_1} \alpha= \int_{D_1}\beta-\int_{L_1}\alpha=
\int_{D_2}\beta-\int_{L_2}\alpha\leq {\rm
Vol}(D_2)-\int_{L_2}\alpha={\cal E}(D_2, L_2) .
\label{stuffc}
\eeq
 The first
equality follows from the definition of the functional. The second
equality follows from the assumption that $X$ is calibrated. The
third equality follows because the class $(\alpha, \beta)$ is
trivial in $H^{k+1}_{dR}(M,N)$. Finally the inequality follows
from the definition of generalized calibration. \rightline{{\bf
Q.E.D.}}

Next we shall investigate some of the properties of the functional
(\ref{mfun}). In particular we have the following:
\begin{theorem}
Let $K$ be a closed submanifold of $N\subset M$ and $L$ a submanifold
of $M$ such that $\partial L=K$. In
addition assume that $(\alpha, \beta)$ represent the trivial
class in $H_{dR}(M,N)$. Then the functional ${\cal E}(K,L)$
is independent of the choice of $L$.
\end{theorem}
{\it Proof:} Let $L'$ another submanifold of $M$ such that
 $\partial L'=K$. Then we have
\begin{equation}
{\cal E}(K, L')={\rm Vol}(K)-\int_{L'} \alpha={\rm Vol}(K)
-\int_K \gamma
={\rm Vol}(K)-\int_{L} \alpha= {\cal E}(K,L)\ .
\end{equation}
The first equality follows from the definition. The second equality
follows from Stoke's theorem because $\alpha$
is exact and so $\alpha=d\gamma$. The third equality also follows
 for the same reason as the second, and the last
follows from the definition of the functional ${\cal E}$.
\hfill{\bf Q.E.D.}

\begin{remark}
Suppose that $(\alpha,\beta)$ represents a class in
 $H^{k+1}(M,N;\bZ)$. Then
it is straightforward to see that
the functional ${\cal E}(K,L) \ \rm {mod} \ \bZ$ is independent of
the choice of $L$.
\end{remark}

\section{Special Cases}

A special case of interest is whenever the generalized calibration
 form $\beta$ is defined as a generalized calibration on $M$.
 For $(\alpha, \beta)$
 to be a trivial class, $\alpha=d\gamma$, $\alpha=d\beta$ and
 $\gamma-\beta=d\zeta$ must be
 an exact form in $M$. In such a case the functional ${\cal E}$
 can be written as
 \begin{equation}
{\cal E}(D)={\rm Vol}(D)-\int_D \gamma\ .
 \end{equation}
 If $X$ is a  compact calibrated submanifold of $M$ with
 boundary $\partial X$, then
 \begin{equation}
{\cal E}(X)=-\int_{\partial X} \zeta\ .
 \end{equation}
In particular if $X$ is closed, then ${\cal E}(X)=0$.
 Since ${\cal E}(Z)>0$ if $Z$
is a closed but not calibrated submanifold of
$M$, calibrated closed submanifolds of $M$ are  global
 minima of ${\cal E}$.
Another special case to consider is whenever $\alpha=d\gamma$,
and we
choose $\gamma=\beta$. This is the case which
we shall focus on later. The functional ${\cal E}$ in this case is
\begin{equation}\label{kfun}
 E(D)={\rm Vol}(D)-\int_D \beta\ .
 \end{equation}
It is worth adapting the main theorem of
generalized calibrations to this case. In particular we have the
following:

\begin{theorem}

Let $\beta$ be a generalized calibration in $M$; then
calibrated submanifolds $X$ of $M$ minimize the functional $E$ in
(\ref{kfun}).

\end{theorem}

{\it Proof:} Let $X$ be a calibrated submanifold of $M$ and
$D$ an open ball in $X$. Next let $D'$ an open ball
in $M$, such that $\partial D=\partial D'$. Then we have
\begin{equation}
E(D)={\rm Vol}(D)-\int_D \beta=\int_D \beta-\int_D \beta
=0=\int_{D'} \beta-\int_{D'} \beta\leq {\rm
Vol}(D')-\int_{D'} \beta= E(D')\ .
\end{equation}
The equalities are obvious. The inequality follows from the
defining property of the calibration form.
\hfill{{\bf Q.E.D.}}

 The energy  (\ref{kfun})  vanishes when evaluated at
 every calibrated submanifold $X$ of $M$. In addition $E(Z)>0$
 if $Z$ is not calibrated and therefore calibrated submanifolds
 are global minima of $E$. From now on, we shall
 focus on the calibrated submanifolds which are the minima
 of the functional (\ref{kfun}).

\begin{remark}
 Standard calibrations for
 which the calibration form is closed $d\beta=0$ are special cases
 of the generalized calibrations associated with the functional
 (\ref{kfun}).
 The only difference is that the energy functional used for
 standard calibrations is the induced volume ${\rm Vol}$.
 Calibrated
 submanifolds under the standard calibrations are minimal.
 For the generalized calibrations the functional (\ref{kfun})
 is not the induced volume but nevertheless it can be identified
 with the
 \lq\lq energy" of the submanifold.
 Observe that the relation between (\ref{kfun}) and induced
 volume evaluated on a closed submanifold $Z$ of $M$ is
 \begin{equation}
 E(Z)={\rm Vol}(Z)-\beta[Z]\ .
 \end{equation}
 The last term depends only on
 the cohomology class of $\beta$. In particular it
 does not contribute in the equations for the criticality of
 $E$ and so the generalized calibrated submanifolds $X$ are minimal.
Conversely, if $\beta$ is a closed form, then minimal submanifolds
of $M$  minimize the energy functional
 $E$.
 \end{remark}

\section{Deformation of generalized calibrations}

\subsection{The second variation of the energy functional}

Here we compute  the second variation of the energy functional
evaluated on a calibrated submanifold. Suppose
 $X$ is a calibrated submanifold of $M$ with respect to a
 generalized calibration $\f$,
 and let $X(t)$ be a 1-parameter
family of submanifolds of $M$ with $X(0)=X$.  In addition let
$d{\rm vol}(t) \equiv d{\rm vol}(X(t))$ be the
volume form of $X(t)$, and $\f(t)$ be the restriction of the
calibration form on $X(t)$. Since $d{\rm vol}(t)$ and
$\f(t)$ are top forms on $X(t)$, we have $\f(t) = \la(t) d{\rm vol}(t)$
for some function $\la(t)$ where $\la(0)=1$.
 This condition follows because the volume form of a calibrated
submanifold is equal to the calibration form. Setting $E(t)=E(X(t))$,
we have  that the  energy
 functional (\ref{kfun}) is
 \bea
 E(0)&=&\int_{X(0)} \phi-\int_{X(0)} \phi
\nonumber \\
&=&\int_{X(t)} \f- \int_{X(t)} \f
\nonumber \\
&=&\int_{X} \f(t) - \int_{X} \f(t)
\nonumber \\
&=&\int_X\big( \la(t) d{\rm vol}(t) - \f(t)\big) \ .
\label{identa}
 \eea
Thus $ \int_X \la(t) d{\rm vol}(t) - \f(t)$ is independent of $t$.
Differentiating the energy functional twice and evaluating at $t=0$,
 we obtain
\begin{equation}
 {d^2 \over dt^2}E(t)|_{t=0}=  -\int_X
{d^2\over dt^2}\la(t)|_{t=0} d{\rm vol}(0)
 -2 \int_X {d\over dt}\la(t)|_{t=0} {d\over dt} d{\rm vol}(t)|_{t=0} \ ,
\label{identc}
\end{equation}
where we have used $\la(0)=1$ and (\ref{identa}). To
proceed with the computation of the
second variation of the energy functional, we shall show
that ${d\over dt}\la(t)|_{t=0}=0$. This a consequence
of the calibration bound. In particular we have the following:
\begin{equation}
d{\rm vol}(t)-\f(t)= (1-\la(t)) d{\rm vol}(t)\ .
\end{equation}
Evaluating this on an appropriately oriented orthonormal basis and
using the calibration bound,
we find that $\la(t)\leq 1$. Since $\la(0)=1$ is a maximum, we have
 ${d\over dt}\la(t)|_{t=0}=0$.

To conclude, the second variation of the energy is
\begin{equation}
 {d^2 \over dt^2}E(t)|_{t=0}=  -\int_X
{d^2\over dt^2}\la(t)|_{t=0} d{\rm vol}(0)  \ .
\label{identd}
\end{equation}

As in the case of standard calibrations, the second variation of the
energy can be computed in terms of the normal vector field $V$.
 The proof
is similar to that given in \cite{ML}. The result is summarized in the
following theorem:
\begin{theorem}
\begin{equation}
 {d^2 \over dt^2}E(t)|_{t=0}=\int_X \big(||\nabla^{\perp} V||^2
d{\rm vol}(0)
 -\nabla^{g}_V\nabla^{g}_V\phi
 -i_{\nabla^{\perp} V}i_{\nabla^{\perp} V}\phi
 -2i_{\nabla^{\perp} V}\nabla^{g}_V\phi\big)
\label{seccvar}
 \end{equation}
 where $\nabla^{g}_V$ is the Levi-Civita covariant derivative of
 $M$ along the
 normal direction $V$ of $X$ and $\nabla^{\perp} V$ is the
 covariant derivative
 of the normal bundle of $X$ in $M$ induced by the Levi-Civita
connection of $M$.
\end{theorem}

\subsection{Special holonomy and generalized calibrations}

Let $(M,g)$ be a Riemannian manifold which admits a metric
connection $\nabla$ with possibly non-vanishing torsion and holonomy
contained in one of the groups $U(n)$ ($2n$), $SU(n)$ ($2n$),
 $Sp(n)$ $(4n)$,
 $Sp(1)\cdot Sp(n)$ $(4n)$, $G_2$ $(7)$ and
$Spin(7)$ $(8)$; the entry in $(\cdot)$ is the real dimension of $M$.
Manifolds with such holonomy admit
generalized calibration forms. These are forms parallel with
respect to the connection $\nabla$. In what follows we
 shall not investigate all cases. Rather we shall focus on
Riemannian manifolds $(M,g)$ which admit a metric connection
$\nabla$ with possibly non-vanishing torsion and holonomy $U(n)$
($2n$), $SU(n)$ ($2n$), $G_2$ $(7)$ and $Spin(7)$ $(8)$. The
general theory of deformations of generalized calibrations will be
developed without further assumptions. However
 in many examples that we shall
present later, we shall require that $(M,g)$ satisfies some additional
geometric conditions in addition to those that are a consequence
of the reduction of the structure group of $TM$.
These will simplify some
aspects of the deformation theory of generalized calibrations
and in particular
the deformation equations.  In particular we shall consider the
following cases:

\vskip 0.4cm
\leftline{{\underline {Holonomy $U(n)$}}}
\vskip 0.3cm

Suppose that a Riemannian  manifold $(M,g)$ (${\rm dim} \ M=2n$) is
equipped with a metric connection $\nabla$ whose holonomy is
contained in $U(n)$. Then $M$ admits an {\it almost} complex
structure $J$, $J^2=-1$, which is parallel with respect to
$\nabla$ and the metric $g$ is hermitian with respect to $J$,
$g(JX,JY)=g(X,Y)$ for $X,Y$ vector fields on $M$. Therefore
$(M,g,J)$ is an {\it almost hermitian} manifold with compatible
connection $\nabla$. Conversely, let $(M,J, g)$ be an almost
hermitian manifold, then $(M,J, g)$ admits a connection $\nabla$
with holonomy contained in $U(n)$. Such a connection $\nabla$ can
be constructed from the Levi-Civita connection $\nabla^{g}$ of $g$
as
\beq \nabla_XY=\nabla_X^{g}Y-{1\over2}J(\nabla_X^{g}J)Y\ ,
\label{ahc} \eeq
 where $X,Y$ are vector fields on $M$. It is
straightforward to verify that $g$ and $J$ are $\nabla$-parallel,
$\nabla g=\nabla J=0$. Note that any almost hermitian manifold $(M,g,J)$
has a K\"ahler form $\Omega(X,Y)=g(X, JY)$. $\Omega$ is
 $\nabla$-parallel
but it is not closed, $d\Omega\not=0$.

\begin{proposition}

Let $(M,g,J)$ be an almost hermitian manifold and $\Omega$
 be the associated
K{\"a}hler form. The forms
$\phi_k={1\over k!} \Omega^k$ are generalized calibrations of
degree 2k. The contact set at every point of $M$ is $Gr(k, \bC^n)$.

\end{proposition}

{\it Proof:} To show this, we shall demonstrate the above
statement at a neighbourhood $U$ of a point $p\in M$. Then because
$\phi_k$ is parallel, it will hold everywhere in $M$. We remark
that there is a neighbourhood of a point $p\in M$ and a local
frame $\{e^a, \e^{\bar a}; a=1,\dots,n\}$, ($e^{\bar
a}=\bar{(e^a)}$), of $(M,g,J)$ such that the metric and K{\"a}hler
form can be written as \bea g&=&\sum^n_{a,\bar b=1}\delta_{a\bar
b} e^a e^{\bar b} \cr \Omega&=&-i \sum^n_{a,\bar b=1}\delta_{a\bar
b} e^a\land e^{\bar b}\ . \eea Observe that in such a frame the
metric and K{\"a}hler form  take the standard form of a Euclidean
metric and (almost) complex structure on $\bR^{2n}=\bC^n$. It
follows that $\phi_k$ are calibrations from  Wirtinger's
inequality on $\bR^{2n}$. From the same inequality it follows that
the contact set consists of the complex k-planes in $\bR^{2n}=\bC^n$.
 All the planes of the contact set can be constructed by
acting with $U(n)$ on the k-plane
\beq
 \xi_0=\{(z_1,\dots, z_k, 0,\dots,
0):z_1,\dots,z_k\in \bC\}\subset \bC^n \ .
\label{stuffe}
\eeq \rightline{{\bf Q.E.D.}}

The calibrated submanifolds $X$ are almost hermitian submanifolds
of $(M,g,J)$. Both the metric and almost complex structure on $X$
are induced from those on $M$; the almost complex structure on $X$
is induced from that on $M$ because at every point $p\in X$, the
holomorphic subspace of $T_pX\otimes \bC$ is identified with a
complex k-plane of the contact set at $p$. Observe that the
dimension of the contact set at every point is $2k (n-1)$. We
shall refer to these generalized calibrations as {\it almost
hermitian calibrations}. These results can be summarized as
 follows:

\begin{corollary}
The almost hermitian calibrations of degree $2k$  of an almost
hermitian manifold $(M,g,J)$ are almost hermitian
 submanifolds of real dimension $2k$.
 \end{corollary}

Suppose that $(M,g,J)$ is a {\it hermitian} manifold, i.e. the
almost complex structure $J$ is integrable.  It is known that such
manifolds admit various connections $\nabla$ with non-vanishing
torsion such that $\nabla g=\nabla J=0$. Because of this the
holonomy of all such connections $\nabla$ is contained in $U(n)$.
Again the forms $\phi_k={1\over k!} \Omega^k$ are generalized
calibrations. The contact set at every point of $M$ is $Gr(k,
\bC^n)$. The proof is identical to the one given above for the
almost hermitian manifolds. In this case, the calibrated
submanifolds $X$ are hermitian submanifolds of $(M,g,J)$. In
particular, they are complex submanifolds. To show this observe
that the Nijenhuis tensor of $X$ vanishes because the complex
structure $J$ of $M$ is integrable. Such calibrations have been
called {\it hermitian calibrations} in \cite{gutpap}. These
results can be summarized as
 follows:

 \begin{corollary}
The  hermitian calibrations  of degree $2k$ of a  hermitian
manifold $(M,g,J)$ are hermitian
 submanifolds of real dimension $2k$.
 \end{corollary}

\vskip 0.4cm
\leftline{{\underline {Holonomy $SU(n)$}}}
\vskip 0.3cm

Suppose that a Riemannian  manifold $(M,g)$ (${\rm dim} \ M=2n$) admits
a metric connection $\nabla$ whose holonomy
is contained in $SU(n)$. In such a case $(M,g)$ admits
an almost complex structure $J$ such that
  $(M,g,J)$ is an almost hermitian manifold  equipped
  with a (n,0)-form $\psi$ such that
$\nabla g=\nabla J=\nabla\psi =0$.  Compatibility of these
conditions requires that the form $\psi$
 is appropriately normalized. In particular
\beq
 (-1)^{{1\over2} (n-1) n} ({i\over2})^n \psi\wedge \bar \psi=
  d{\rm vol}\ ,
 \label{nvol}
 \eeq
 where
$\bar \psi$ is the associated $(0,n)$ form and $d{\rm vol}$ is the
volume form of
 $M$ with respect to the metric $g$.

Conversely, let $(M,g,J)$ be an almost hermitian manifold with
topologically trivial canonical bundle, then there is a connection
$\nabla$ which has holonomy contained in $SU(n)$. To see this,
observe that  $(M,g,J)$  admits a no-where vanishing (n,0)-form
$\psi$ but in general $\psi\wedge \bar\psi$ will not be related to
the volume form of $M$  as in (\ref{nvol}). In
 general, we have
 \beq
 \psi\wedge \bar \psi= (-1)^{-{1\over2} (n-1) n+n} (2i)^n f^2
 d{\rm vol_g}
 \label{nnvol}
 \eeq
for some nowhere vanishing real-valued function $f$ of $M$, where
$\bar \psi$ is the associated $(0,n)$ form. Now there are two
possibilities to consider. First, define $\chi= f^{-1}\psi$.
Observe that $\chi$ is again a nowhere vanishing section of the
canonical bundle and it is normalized
 as in (\ref{nvol}). In such case, one can show that the connection
\beq
\nabla_iY^j=\nabla^{g,J}_iY^j+{1\over n!}\chi^{jk_1\dots
k_{n-1}}{}\nabla^{g,J}_i\chi_{kk_1\dots k_{n-1}} Y^k+ {1\over
n!}\bar\chi^{jk_1\dots k_{n-1}}{}\nabla^{g,J}_i\bar\chi_{kk_1\dots
k_{n-1}} Y^k
\label{stuffh}
\eeq
has holonomy contained in $SU(n)$, i.e. $g$, $J$ and
$\chi$ are all $\nabla$-parallel, where $\nabla^{g,J}$ is a
connection with holonomy contained in $U(n)$,
$\nabla^{g,J}g=\nabla^{g,J}J=0$,
 such as the one given in (\ref{ahc}).

Alternatively, observe that $f$ is either a positive or negative
function. So without loss of generality we can take $f$ to be
positive because if it is negative we can take $|f|$. Next define
a new metric $h$ on $M$ by $h=f^{2\over n} g$. Then observe that
$(M,h,J)$ is again a hermitian manifold and the (n,0)-form $\psi$
is normalized as in (\ref{nvol}) with respect to the new metric
$h$. In such a case, one can show that the connection
\beq
\nabla_iY^j=\nabla^{h,J}_iY^j+{1\over n!}\psi^{jk_1\dots
k_{n-1}}{}\nabla^{h,J}_i\psi_{kk_1\dots k_{n-1}} Y^k+ {1\over
n!}\bar\psi^{jk_1\dots k_{n-1}}{}\nabla^{h,J}_i\bar\psi_{kk_1\dots
k_{n-1}} Y^k
\label{stuffi}
\eeq
 has holonomy contained in $SU(n)$, ie $h$, $J$ and
$\psi$ are all $\nabla$-parallel, where $\nabla^{h,J}$ is a
connection with holonomy contained in $U(n)$,
$\nabla^{h,J}h=\nabla^{h,J}J=0$, such as the one given in (\ref{ahc})
but constructed using the  Levi-Civita connection of the metric
$h$.

\begin{proposition}
Let $(M,g,J)$ be an almost hermitian manifold, ${\rm dim} M=2n$,
 with trivial canonical
 bundle and associated parallel
(n,0)-form $\psi$. The form $\phi={\rm
Re}(\psi)$ is a generalized calibration of
degree n. The contact set $C_p$ at every point of $(M,g,J)$
is $SU(n)/SO(n)$.
\end{proposition}

{\it Proof:} As in the case of almost hermitian calibrations, we
shall prove the proposition in a neighbourhood $U$ of a point $p\in
M$.  In such a neighbourhood, there is a local frame $\{e^a,
e^{a'}\}$ in $(M,g,J)$ such that \bea g&=&\sum_{a=1}^n
\big((e^a)^2+(e^{a'})^2\big) \cr \Omega&=& \sum^n_{a=1} e^a\land
e^{a'} \cr
 \psi&=& (e^1+i e^{1'})\land \dots \land (e^n+i e^{n'})\ .
\eea Again in this frame the metric, K{\"a}hler form and
(n,0)-form take the standard form of those in $\bR^{2n}$. These
are precisely the data of a SLAG calibration in $\bR^{2n}$.
Therefore it follows that $\phi$ is a generalized calibration
from the results of Harvey and Lawson as they apply for Special
Lagrangian (SLAG) calibrations. Similarly, it follows that the
contact set at every point $p\in M$ is $SU(n)/SO(n)$. All the
calibrated planes can be constructed by acting with $SU(n)$ on the
standard plane
\beq
 \xi_0=\{(x_1, \dots, x_n, 0,\dots,0): x_1,
\dots, x_n\in \bR\}\subset \bR^{2n}\ .
\label{stuffj}
\eeq \rightline{{\bf Q.E.D.}}

Such generalized calibrations have been called {\it Special Almost
Symplectic} or SAS for short \cite{gutpap}. The SAS calibrations
are real middle dimension submanifolds of $(M,g,J)$. We remark
that if ${\rm Re}\psi$ is a SAS calibration, then ${\rm Re
}(e^{i\theta}\psi)$ is also a SAS calibration, where $\theta$ is a
constant angle.

So if $(M,g,J)$ is an almost hermitian manifold with a
compatible connection which
has holonomy contained in $SU(n)$, then $M$ admits two
 types of generalized calibrations
with calibrated submanifolds;
the almost
 hermitian and the SAS. This is reminiscent
of Calabi-Yau manifolds which have two types of
calibrated submanifolds the
K{\"a}hler and SLAG.

A special case that we shall investigate later is that in
which $(M,g,J)$ is a hermitian manifold which admits a compatible
 connection
$\nabla$ with holonomy contained in $SU(n)$.
Again these manifolds admit two types of calibrations; the hermitian and
the SAS. The contact set of SAS calibrations is $SU(n)/SO(n)$ at
 every point of $M$.

The above results are summarized as follows:

\begin{corollary}
Let $(M,g,J)$ be an (almost) hermitian manifold with trivial
 canonical bundle.
Then $(M,g,J)$ admits (almost) hermitian and SAS calibrations.
\end{corollary}

Another class of hermitian manifolds $(M,g,J)$ are those for which
the canonical bundle is holomorphically trivial. Connections on
such manifolds will be investigated in the second part of this
paper. Such manifolds admit a holomorphic $(n,0)$-form $\psi$.
Thus $\psi$ is closed, $d\psi=0$, but $\psi$ is not always
normalized as in (\ref{nvol}). In such a case, as we have
explained, we can either rescale the form $\psi$ or we can rescale
the metric $g$. In particular  we can rescale the metric $g$ as
$h=f^{2\over n} g$ so that $\psi$ becomes a calibration form,
where $f$ is given in (\ref{nnvol}). Since $\psi$ remains closed,
the associated calibrated submanifolds are SLAGs and therefore
minimal with respect to $h$. Alternatively as we have seen, one
can rescale the holomorphic $(n,0)$-form $\psi$ as
$\chi=f^{-1}\psi$ where $f$ is given in (\ref{nnvol}). In this
case, the rescaled form $\chi$ is a calibration but it is not
closed.  The associated calibrated submanifolds are SAS with
respect to the original metric $g$.

Another special case that has recently been investigated  is
that of K\"ahler manifolds $(M,g,J)$ with trivial canonical bundle
for which the metric $g$ is not a Calabi-Yau metric. Such
manifolds have been called {\it almost Calabi-Yau} and have been
studied
in the context of mirror symmetry (see \cite{Jo1}). Note that as a
consequence of the Calabi-Yau theorem, compact almost Calabi-Yau
manifolds always admit a Calabi-Yau metric. Compact almost
Calabi-Yau manifolds admit a holomorphic $(n,0)$ form $\psi$,
which is therefore closed ($d\psi=0$), but not necessarily parallel
 with
respect to a hermitian connection because it does not satisfy the
normalization condition (\ref{nvol}). The strategy adopted in this
case is to conformally rescale the metric $g$, as in the case of
hermitian manifolds with holomorphically trivial canonical bundle
above, so that $\psi$ remains closed, and the associated
calibrated submanifolds are Special Lagrangian and therefore
minimal with respect to the rescaled metric. Alternatively as we
have seen, one can rescale the holomorphic $(n,0)$-form $\psi$.
The associated calibrated submanifolds are SAS with respect to the
original metric $g$.

In the case that $(M,g,J)$ is Calabi-Yau, then the hermitian
 calibrations
become the standard K{\"a}hler calibrations while the SAS
 calibrations
become the standard  SLAG calibrations.

\vskip 0.4cm
\leftline{{\underline {Holonomy $G_2$}}}
\vskip 0.3cm

A Riemannian manifold $(M,g)$  (${\rm dim} \ M=7$) equipped
with a metric connection $\nabla$ whose holonomy is
contained in $G_2$ admits a
 $\nabla$-parallel three-form $\psi$ and a $\nabla$-parallel
 four-form
$*\psi$ which is the dual of $\psi$. As we shall show these forms
$\psi$ and $*\psi$  are  generalized associative and generalized
co-associative calibrations, respectively.

Conversely, let $(M,g,\psi)$ be a Riemannian manifold which admits
a three-form $\psi$ that satisfies the algebraic conditions of a
$G_2$ invariant structure (such a three-form is stable in the
terminology of \cite{Hit}),
 then there is a connection $\nabla$ which has holonomy $G_2$.
This connection can be expressed in terms of the Levi-Civita
connection $\nabla^g$ of $g$ and the form $\psi$ as
\beq
 \nabla_k
Y^i=\nabla_k^g Y^i+ {1\over 18} \psi^i{}_{pq} \nabla^g_k
\psi_j{}^{pq} Y^j+ {1\over 108} *\psi^i{}_{pqr} \nabla^g_k
*\psi_j{}^{pqr} Y^j\ ,
\label{stuffk}
\eeq
where $Y$ is a vector field.

\begin{proposition}
Let $(M,g,\psi)$ a seven-dimensional manifold which admits
a $G_2$-structure as above. The forms
$\psi$ and $*\psi$ are generalized calibrations of
degree three and four, respectively. In both cases the
contact set at every point of $M$ is $G_2/SO(4)$.
\end{proposition}

{\it Proof:} To show this, we remark that there is locally a
frame $\{e^A; A=1,\dots, 7\}$ of $(M,g,J)$ such that
the metric and the $\nabla$-parallel three-form $\psi$
 can be written as
\bea
g&=&\sum_{A=1}^7 (e^A)^2
\cr
\psi&=& e^{123} + e^1 \land (e^{45}-e^{67}) +e^2 \land
(e^{46}+e^{57}) +e^3 \land (e^{47}-e^{56})\ ,
\label{gstring}
\eea
where $e^{12}=e^1 \land e^2$ and similarly for the rest.
Observe that in such a frame the metric $g$ and
parallel three-form $\psi$ take  the standard form of a
Euclidean metric and $G_2$-invariant three-form in $\bR^7$.
It follows that both $\psi$ and $*\psi$
are generalized calibrations from the results
of Harvey and Lawson as they apply for
associative and co-associative calibrations. It also follows
that the contact set at every point of $M$ is $G_2/SO(4)$.
All the planes of the contact set
can be constructed by acting with $G_2$
on the 3-plane
\beq
\xi_0=\{(x_1,x_2, x_3, 0, 0, 0, 0):x_1, x_2, x_3\in \bR\}
\subset\bR^7  \ ,
\label{stuffl}
\eeq
for generalized associative calibrations and on the 4-plane
\beq
\xi_0=\{( 0,0, 0, x_4,x_5, x_6, x_7): x_4,x_5, x_6, x_7\in
 \bR\}\subset \bR^7
 \ ,
\label{stuffm}
\eeq
for generalized co-associative calibrations.
\rightline{{\bf Q.E.D.}}

There are many special cases of seven-dimensional Riemannian
manifolds that admit connections whose holonomy is contained in
$G_2$. We shall present many of these cases when we investigate
the deformation theory of generalized associative and
co-associative calibrations.

\vskip 0.4cm
\leftline{{\underline {Holonomy $Spin(7)$}}}
\vskip 0.3cm

A Riemannian manifold $(M,g)$  (${\rm dim} \ M=8$) equipped with a
metric connection $\nabla$ whose holonomy is contained in
$Spin(7)$ admits a $\nabla$-parallel self-dual four-form $\Phi$.
As we shall see, $\Phi$ is a generalized Cayley calibration.

Conversely, let $(M,g,\Phi)$ be a Riemannian manifold equipped
with a self-dual four-form $\Phi$ which satisfies the algebraic
conditions of a $Spin(7)$ structure, then $(M,g,\Phi)$ admits a
connection $\nabla$ whose holonomy is contained in $Spin(7)$. The
connection $\nabla$ can be expressed in terms of the Levi-Civita
connection $\nabla^g$ of $g$ as
\beq \nabla_k Y^i=\nabla_k^g
Y^i+{1\over 96}  \Phi^i{}_{mkl} \nabla^g_k \Phi_j{}^{mkl} Y^j\ ,
\label{stuffn}
\eeq
 where $Y$ is a vector field.

\begin{proposition}
Let $(M,g,\Phi)$ be an eight-dimensional manifold with $Spin(7)$
 structure as above.
 The forms $\Phi$ is generalized calibration of degree four. In
both cases the contact set at every point of $M$ is $Spin(7)/K$,
where $K=SU(2)\times SU(2)\times SU(2)/Z_2$.
\end{proposition}

{\it Proof:}  To show this, we remark that there is locally a frame
$\{e^A; A=1,\dots, 8\}$ of $(M,g,J)$ such that the metric and the
$\nabla$-parallel self-dual four-form $\Phi$ can be written as
\bea
g&=&\sum_{A=1}^8 (e^A)^2 \cr \Phi &=& e^{1234} +
(e^{12}-e^{34}) \land (e^{56}-e^{78})
 +(e^{13}+e^{24}) \land (e^{57}+e^{68})
\cr
&+&(e^{14}-e^{23}) \land (e^{58}-e^{67}) +e^{5678}\ .
\label{sdff}
\eea Observe that in such a frame the metric and parallel
self-dual four-form take  the standard form of a Euclidean metric
and $Spin(7)$-invariant four-form in $\bR^8$. It follows that
$\Phi$ is a calibration from the results of Harvey and Lawson as
they apply for the Cayley calibration. It also follows that the
contact set at every point of $M$ is $Spin(7)/K$, where
$K=SU(2)\times SU(2)\times SU(2)/Z_2$. All the planes of the
contact set can be constructed by acting with $Spin(7)$ on the
4-plane
\beq
 \xi_0=\{(x_1,x_2, x_3, x_4, 0,\dots, 0): x_1,x_2, x_3,
x_4\in \bR\}\subset \bR^8 \ .
\label{stuffo}
\eeq \rightline{{\bf Q.E.D.}}

It can be shown that given a self-dual four-form as in
(\ref{sdff}), there is always a connection $\nabla$ with torsion a
three-form given in \cite{siv} which has holonomy contained in
$Spin(7)$. The torsion is
\beq
 T=\delta\Phi+{7\over6} *(\theta\land
\Phi)
\label{stuffp}
\eeq
 where $\theta={1\over7} *(\delta \Phi\land \Phi)$ is the
Lee form of the manifold.

\subsection{Useful Formulae}

For the investigation of the deformation
theory of generalized calibrations,
we shall use some formulae which relate the
Lie derivative of a form to a covariant
derivative.

 Let $\chi$ be a k-form
expressed as
 $\chi={1\over k!} \chi_{A_1\dots A_k} e^{A_1}\land \cdots \land
 e^{A_k}$
 in a frame $\{e^{A_p}\}$. Then
 \beq
 d\chi ={1\over k!} \nabla_{A_1}
 \chi_{A_2 \dots A_{k+1}} e^{A_1}\land \cdots \land e^{A_{k+1}}
 +{1\over (k-1)!} \chi_{A_1\dots A_k} T^{A_1}\land e^{A_2}\land
\cdots \land e^{A_k}
 \label{excon}
 \eeq
where $T^A=\nabla e^A={1\over2} T^A{}_{BC} e^B\land e^C$ is the
torsion 2-form of $\nabla$ or equivalently
\beq
 T(X,Y)= \nabla_X Y -
\nabla_Y X - \big[ X \ , \ Y \big]
\label{stuffq}
\eeq
 for vector fields $X$, $Y$.

The Lie derivative of $\chi$ with respect to a vector field
$V=V^A e_A$,  $(e^A, e_B)=\delta^A{}_B$,
 is as follows:
  \bea
{\cal L}_V\chi&=&(i_Vd+di_V)\chi={k+1\over k!} V^{A_1}
 (\nabla\chi)_{A_1A_2\dots A_{k+1}}
 e^{A_2}\land \cdots \land e^{A_{k+1}}\nonumber\\
 &+&{1\over (k-1)!} \chi_{A_1A_2\dots A_k} i_V T^{A_1}
\land e^{A_2}\land
  \cdots \land e^{A_k}\nonumber\\
& +&{1\over (k-2)!} V^{B}\chi_{A_1 B A_3\dots A_k}  T^{A_1}
\land e^{A_3}\land
  \cdots \land e^{A_k}\nonumber\\
&+&{1\over (k-1)!} \nabla_{A_1}(V^{B} \chi_{B A_2\dots A_k})
 e^{A_1}\land e^{A_2}\land \dots\land e^{A_k}\nonumber\\
& +&{1\over (k-2)!} V^B \chi_{BA_2 A_3\dots A_k}
T^{A_2}\land e^{A_3}\land\dots\land e^{A_k}\nonumber\\
&=& {1\over k!} V^B \nabla_B\chi_{A_1\dots A_k} e^{A_1}
\land\dots\land e^{A_k}\nonumber\\
&+&{1\over (k-1)!} \tilde\nabla_{A_1}
 V^B \chi_{BA_2\dots A_k} e^{A_1}\land \dots\land e^{A_k}\ ,
 \eea
 where $\tilde\nabla$ is the unique connection associated
 with $\nabla$ which
 has torsion $\tilde T=-T$. To summarize
 \beq
 \cL_V \chi= \nabla_V \chi +{1 \over (k-1)!}
\chi_{A A_1 \dots A_{k-1}}
\tilde \nabla_B V^A e^B \land
 e^{A_1} \land \dots \land e^{A_{k-1}}\ .
 \label{ppara}
 \eeq
  So, if $\chi$ is $\nabla$-parallel, $\nabla\chi=0$, and
 \beq
 \cL_V \chi= {1 \over (k-1)!} \chi_{A A_1 \dots A_{k-1}}
 \tilde \nabla{}_B V^A e^B \land
 e^{A_1} \land \dots \land e^{A_{k-1}} \ .
 \label{lie}
 \eeq

Another formula which we shall find useful is the Lie derivative of
a vector-valued k-form
\beq
\xi = {1 \over k^!} \xi_{A_1 \dots A_k}{}^B e^{A_1}
\land \dots \land
e^{A_k} \otimes e_B
\label{stuffr}
\eeq
which may be written as
\beq
\cL_V \xi = \nabla_V \xi +{1 \over (k-1)!} \xi_{B A_2 \dots A_k}{}^A
\tilde\nabla_{A_1} V^{B} e^{A_1} \land \dots \land e^{A_k}
\otimes e_A
-{1 \over k!} \xi_{A_1 \dots A_k}{}^B \tilde\nabla_B V^A
 e^{A_1} \land \dots \land e^{A_k} \otimes e_A\ .
\label{stuffs}
\eeq
So, if $\xi$ is $\nabla$-parallel, and $X$ is some
submanifold on which
\beq
\xi_{A_1 \dots A_k}{}^A
 e^{A_1} \land \dots \land e^{A_k}|_X\otimes e_A=0
\label{stufft}
\eeq
then
\beq
\cL_V \xi |_X = {1 \over (k-1)!} \xi_{B A_2 \dots A_k}{}^A
\tilde\nabla_{A_1} V^{B} e^{A_1} \land \dots \land e^{A_k}|_X
 \otimes e_A
\ .
\label{lievf}
\eeq
These formulae are key in the investigation of the
deformation theory of generalized calibrations.

\subsection{Second variation of the energy functional revisited}

The second variation of the energy functional is considerably
simplified if we assume that the calibration form $\f$ is
$\nabla$-parallel.  In particular, suppose that $X$ is a
k-dimensional sub-manifold calibrated with respect to $\f$.
Suppose we consider an adapted frame so that $\{ e^a : a=1, \dots
,k \} $ are tangent to $X$ and $\{ e^i : i=k+1 , \dots , {\rm{ dim
(M)}} \}$ are normal to $X$, and we take  $\f_{i a_1 \dots
a_{k-1}}|_{t=0}e^i\wedge e^{a_1}\wedge\dots\wedge e^{a_{k-1}}=0$.
Then we deform the calibration and write the calibration form as
\bea \phi&=&\lambda e^1\wedge\dots\wedge e^k+ {1\over (k-1)!}
\f_{i a_1 \dots a_{k-1}}e^i\wedge e^{a_1}\wedge\dots\wedge
e^{a_{k-1}}\\ &+& {1\over (k-2)!} \f_{i ja_1 \dots
a_{k-2}}e^i\wedge e^j\wedge e^{a_1}\wedge\dots\wedge e^{a_{k-2}}+
{\cal O}((e^i)^3) \ .
\label{exex} \eea
Using the fact that $\phi$ is
$\nabla$-parallel, we can compute  ${d^2 \lambda \over d t^2}$ by
acting twice on (\ref{exex}) with $\nabla_V$ where
$V=\partial/\partial t$ is a normal vector field. After some
lengthy computation, we find \bea {d^2 \over dt^2}E(t)|_{t=0} &=&
-\int_X {d^2\over dt^2}\la(t)|_{t=0} d{\rm vol}(0) \nonumber
\\
&=& \int_X \big[ (\tilde\nabla^\perp V, \tilde\nabla^\perp V)
 d{{\rm vol}}(0)-
i_{\tilde\nabla^\perp V } i_{\tilde\nabla^\perp V}\f \big]|_{t=0}\ .
\label{ssecb}
\eea

\section{Deformation of hermitian calibrations}

Let $(M,g,J)$ be a hermitian manifold of complex dimension $n$.
As we have mentioned the calibrated submanifolds
with respect to $\phi={1\over k!}\Omega^k$ are the complex
 submanifolds of
$(M,g,J)$ of dimension $k$. The deformation
theory of a complex submanifold $X$ of $M$ is well known. The
dimension of the moduli space is the number of
holomorphic vector fields of the normal bundle of the
submanifold $X$, i.e. it is the dimension of the
{\v C}ech cohomology  \v H${}^0(N_X)$ where $N_X$ is the normal
 bundle of $X$ in $M$. There is an obstruction of integrating
these small deformations. This lies in the
{\v C}ech cohomology group \v H${}^1(N_X)$.
It is clear that there should be a theory of deformations
of almost hermitian calibrations. This will be investigated elsewhere.

\subsection{Complex submanifolds of Hermitian manifolds}

It has been known for some time that hermitian manifolds admit compact
complex submanifolds which
represent the trivial homology class. This is unlike the case of
compact K{\"a}hler manifolds where
complex submanifolds always represent a non-trivial homology
class.
An example of such a hermitian manifold that admits a holomorphic
submanifold which represents the trivial homology class has been
given in \cite{alggeom}. This is an example of a hermitian
calibration which is not K{\"a}hler.

Another example of a hermitian calibration is that of the Hopf
fibre $S^1\times S^1$ in the group manifold $S^3\times S^3$. As we
shall demonstrate later in the investigation of examples of SAS
calibrations in group manifolds, such a submanifold is holomorphic
with respect to a hermitian structure on $S^3\times S^3$. Observe
that the hermitian calibration $S^1\times S^1$ represents the
trivial homology class.

\section{Deformations of SAS Calibrations}

Let $(M,g,J)$ be an almost complex manifold which admits a
compatible connection $\nabla$ $(\nabla g=\nabla J=0)$
which has holonomy contained in $SU(n)$. The following can be
shown using the results of Harvey and Lawson:
\begin{theorem}
Let $(M,g,J)$ be a manifold as above. If $X$ is a SAS calibrated
submanifold of $M$ with respect to $\phi={\rm
Re}~\om$, where $\om$ is the parallel-(n,0) form, then
$\Omega|_X={\rm Im}~\om|_X=0$, where $\Omega$ is the
Hermitian form obtained from $J$. Conversely, if $X$ is a
 middle dimension
 submanifold of $M$ such that $\Omega|_X={\rm
Im}~\om|_X=0$, then $X$ is calibrated with respect to $\phi$.
\end{theorem}

Because of this, the small deformations of $X$ generated
 by sections $V$ of the normal bundle, $N_X$,  of $X$ in
$M$ which preserve the  property that  $X$ is calibrated  satisfy
\begin{equation}
\cL_V \Omega |_X =0\ ,~ \cL_V
{\rm Im}~\om |_X =0 \ . \label{defconda}
\end{equation}

To determine the conditions on $V$ imposed by the above two
conditions, we proceed as follows. We introduce an
orthonormal basis $\{e_a, e_{a'}\}=\{ e_1 , \ \dots \ , e_n , e_{1'} ,
\ \dots \ , e_{n'} \}$ of the tangent
bundle of $M$ and a dual frame $\{e^a, e^{a'}\}
=\{ e^1 , \ \dots \ , e^n , e^{1'} , \ \dots \ , e^{n'} \}$ of
$M$ such
 that the K{\"a}hler form $\Omega$ and the parallel
(n,0)-form take the (canonical) forms
\bea
 \Omega & =& \sum_{a=1, \ b'=1}^n \Omega_{ab'} e^a \land
e^{b'}      =\sum_{a=1, \ b'=1}^n \delta_{ab'} e^a \land
e^{b'}=\sum^n_{a=1}e^a \land e^{a'} \nonumber\\
\om & =&\Pi_{a=1}^n \big(e^a-i J(e^a)\big)= (e^1+ie^{1'})
 \land \dots \land (e^n+ie^{n'})\ .
 \label{flatt}
 \eea
  It is
clear that the non-vanishing components of the
 almost complex structure $J$ in this frame are $J^a{}_{a'}$
 and $J^{a'}{}_a$.

 Restricting the orthonormal  basis $\{e_a, e_{a'}\}$ at a point
 $p$ of the calibrated submanifold $X$, $\{e_a\}$ is a basis
 in the tangent space $T_pX$ and $\{e_{a'}\}$ is a basis in the
  fibre $N_p$ of the  normal bundle $N_X$ of $X$ in $M$. Therefore
 the deformations of $X$ in $M$ are described by vector
 fields $V=V^{a'} e_{a'}$.
 Using this and (\ref{lie}), the conditions (\ref{defconda})
 can be written as follows:
\bea
\tilde\nabla_a V^{b'} \Omega_{b'b} e^a\land e^b&=&0
\nonumber\\
\tilde\nabla_a V^{b'} J^{a}{}_{b'}&=&0\ .
\label{sasc}
\eea
These are viewed as equations for the normal vector field $V$.

There is another way to write the deformation equations
 of SAS calibrations. For this observe that the normal bundle
 $N_X$ and the tangent bundle $TX$ of $X$ are isomorphic,
 $N_X\equiv TX$.
 The isomorphism is induced by the almost complex structure $J$
 as $U=U^a e_a=J(V)=J^a{}_{b'} V^{b'}e_a$.
Using this, the conditions (\ref{defconda}) on the normal
 vector field $V$
can be rewritten as
\begin{equation}
d (i_V \Omega)|_X+ i_V d \Omega |_X =0\ ,
~~d (i_V {\rm Im}~\om )|_X+i_V d {\rm Im}~\om|_X  =0 \ .
\label{defcondb}
\end{equation}

Next using the fact that both $\Omega$ and $\omega$ are parallel with
respect to the connection $\nabla$ and (\ref{excon}), these
two conditions can be expressed as
\beq
D_1U\equiv  dU- U_a (T^a+\hat T^a)=0
\label{sasa}
\eeq
 and
\beq D_0^\dagger U\equiv \delta U +U^a (t_a + \hat t_a)=0 \ ,
\label{sasb}
\eeq
where $\hat T^a= e^b\land e^c \Omega_{bb'} T^{b'}{}_{ca'}
 \Omega^{a'a}$, $t_a= T^b{}_{ab}$ and
$\hat t_a=J^{a'}{}_a T^{b'}{}_{a'b} J^b{}_{b'}$.

Equivalently, the deformation equations (\ref{sasa}) and
(\ref{sasb}) can be written in components as follows: \beq
\partial_{\mu_1} U_{\m_2} -\pd{\m_2} U_{\m_1} - U_\p
 (T^\p_{\m_1\mu_2}+
  \hat T_{\mu_1\mu_2}{}^\p \big)=0
\label{rewrb}
 \eeq
and
 \beq
(\nabla^{g})^\m U_\m- U^\m (t_\m +\hat t_\m)=0 \ .
 \label{rewra}
 \eeq
 We remark that
 the deformations of SAS calibrations in the special case when 
 $M$ is a symplectic manifold, and so $\Omega$ is closed,
 have been considered in \cite{Sal}. It has been shown that
 the moduli are unobstructed, and the dimension of the moduli space
 is $b_1 (X)$. Both expressions (\ref{sasc}) and (\ref{sasa}, \ref{sasb})
 of the deformation equations will be used later in the  examples
 to find the moduli space of SAS calibrations in non-symplectic manifolds.
  
\subsection{SAS calibrations and an elliptic system}

To investigate whether the  differential system  (\ref{sasa}) and
(\ref{sasb})
 has solutions,
consider the following resolution:
\begin{equation}
\Lambda^0(X)\stackrel{D_0}{\to}
\Lambda^1(X)\stackrel{D_1}{\to} \Lambda^2(X)
\end{equation}
where $D_1$ has been defined as above and $D_0$ is the
 adjoint of $D_0^\dagger$.
Clearly we have the adjoint resolution as follows:
\begin{equation}
\Lambda^0(X)\stackrel{D^\dagger_0}{\gets} \Lambda^1(X)
\stackrel{D^\dagger_1}{\gets} \Lambda^2(X) \ ,
\end{equation}
where $D^\dagger_1$ is the adjoint of $D_1$. Next we
 can consider the Laplacian
\begin{equation}
\triangle =D_0
D_0^\dagger+D_1^\dagger D_1\ .
\end{equation}
{}From general elliptic theory, we know that the
 solutions of the equations
 $D_1 U=0$ and $D_0^\dagger U=0$ are zero modes
of the Laplacian $\triangle$.
  Conversely, the zero modes of the Laplacian
$\triangle$ are also solutions
  of the two differential equations. From general
 elliptic theory we also have the following:

  \begin{corollary}
The moduli space of a closed SAS calibration $X$ in $M$, if it
exists, has finite dimension.
\end{corollary}

We shall investigate the elliptic system in more detail
in special cases below.

\section{Special Cases}

There are several different types of almost hermitian structures,
for example given in the Gray-Hervella  classification \cite{GrH}.
We shall not explore all cases here. Instead, we shall focus on
some of these. Some explicit examples will be given later. In what
follows, we shall assume that some  hermitian connections have
holonomy contained in $SU(n)$.

\subsection{Almost Hermitian manifolds with skew torsion}

Let $(M,g,J)$ be an almost hermitian manifold. It has been shown
in \cite{FI} that $(M,g,J)$ admits a unique almost hermitian
connection $\nabla$ with torsion a three-form iff the Nijenhuis
tensor of $J$ is a three-form as well, ie $(M,g,J)$ is a $G_1$
manifold in the Gray-Hervella classification. In that case the
torsion of the connection is \beq T(X,Y,Z)=-d\Omega(JX,JY,
JZ)+N(X,Y,Z) \label{stuffaa} \eeq where
$N(X,Y)=[JX,JY]-[X,Y]-J[JX,Y]-J[X,JY]$ is the Nijenhuis tensor.
 Suppose
 that in addition the holonomy of this connection is contained
in $SU(n)$. In such a case the differential system for SAS
calibrations can be simplified to
\bea dU-U_a (T^a+\hat T^a)&=&0
\cr \delta U+ U^a \hat t_a&=&0\ .
\label{stuffab}
\eea
Despite the
simplification of the second deformation equation, we have not
been able to analyze the system further. However, a special class
of such hermitian manifolds are the Nearly K\"ahler ones. For
these, the deformation equations simplify further.

\subsection{Nearly K\"ahler manifolds}
\vskip 0.3cm

Let $(M,g,J)$ be a Nearly K\"ahler manifold, ie $(M,g,J)$ is an
almost hermitian manifold satisfying $(\nabla^g_XJ)Y+(\nabla^g_YJ)X=0$
, where $X,Y$
are vector fields  on $M$. It is known \cite{Gra} that on a Nearly
K{\"a}hler manifold $(M,g,J)$ the following identities hold:
\beq\label{relnk} 4d\Omega(X,Y,Z)=3N(JX,Y,Z)=-12g((\nabla^g_XJ)Y,Z)
\eeq and \beq\label{nk1} 2g((\nabla^g_U\nabla^g_XJ)Y,Z)=-\big(
g((\nabla^g_UJ)X,(\nabla^g_YJ)JZ)+{\rm cyclic} (X,Y,Z)\big).
\eeq
The Nijenhuis tensor $N$ is a $(3,0)+(0,3)$-form. Nearly K\"ahler
manifolds admit a compatible connection $\nabla$ with torsion a
three-form $T={1\over4}N$ and $\nabla T=0$ \cite{Kir,BM,FI}.

We shall focus  our attention to  six-dimensional nearly K\"ahler
manifolds. Any six-dimensional nearly K\"ahler manifold is
Einstein and of constant type (see \cite{Gra}). This means that
the Ricci tensor, $ Ric^g = \frac{5}{2} a g$ and
\beq\label{nk2}
||(\nabla^g_XJ)Y||^2 = \frac{1}{2} a \cdot
\big(||X||^2 \cdot ||Y||^2 - g^2(X,Y) - g^2(X,JY) \big),
\eeq
where
$a = Scal^g/15$ is a positive constant and $Scal^g$ denotes
the scalar curvature of  $g$. It is clear that the holonomy of the
connection $\nabla$ of any six-dimensional Nearly K\"ahler
manifold which is not K\"ahler is contained in $SU(3)$.

To investigate SAS calibrations in six-dimensional Nearly K\"ahler
manifolds, we shall first prove the following theorem for Lagrangian
submanifolds in Nearly K\"ahler manifolds.

\begin{theorem}

A three-dimensional Lagrangian submanifold $L$ of a
six-dimensional Nearly K\"ahler manifold $(M,J,g)$ is a SAS
calibration and  minimal. Consequently any Lagrangian submanifold
$L$ is orientable.

\end{theorem}

{\it Proof:} To show that any Lagrangian submanifold of a
 six-dimensional
nearly K\"ahler manifold is minimal, we shall first show that
\beq\label{par}
g((\nabla^g_XJ)Y,Z)=0
\eeq
 for $X,Y,Z$ tangent to $L$. To see this, we use the fact
 that $J$ is parallel
 with respect to the connection $\nabla$ with torsion
the Nijenhuis tensor $N$.
 Using the fact that $N$ is a $(3,0)+(0,3)$-form and
 that $L$ is Lagrangian, it
 is straightforward  to verify (\ref{par}).

To show that $L$ is SAS, observe that a $\nabla$-parallel (3,0)
form $\psi$ can be defined on $L$ with ${\rm Re}\psi=N$ and ${\rm
Im}\psi={3\over4} d\Omega$; $\Omega$ is the K\"ahler form. Since
both the Nijenhuis tensor $N$ and the almost complex structure $J$
are $\nabla$-parallel, in view of (\ref{relnk}), $\psi$ is
parallel as well. So $N$ can be identified with the calibration
form. Using (\ref{par}), $(d\Omega)|_L=0$ and so $L$ is a SAS
calibration.  Consequently $L$ is orientable.

The second part of the proof of this theorem that $L$ is minimal
is a generalization of a similar theorem for $S^6$ in \cite{Ej}.
To begin, denote with $\alpha$ and $ A$ the second fundamental
form and the shape operator of the submanifold $L$ in a manifold
$M$, respectively. From the definition of $\alpha$ and $A$, we
have   $X,Y$ tangent to $L$ and $\xi$ normal to $L$
\beq
\nabla^g_XY=\nabla^{g_L}_XY+\alpha(X,Y),\quad
\nabla^g_X\xi=-A_{\xi}X+D_X\xi,
\label{stuffac}
\eeq
where $X,Y$ are vector fields
tangent to $L$, $\xi$ is a vector field normal to $L$,
$\nabla^{g_L}$ is the induced Levi-Civita connection on $L$ and
$D$ is the induced connection on the normal bundle. Recall that
$\alpha(X,Y)=\alpha(Y,X)$  and $g(\alpha(X,Y),\xi)=g(A_{\xi}X,Y)$.

To see that $L$ is minimal we observe that the (3,0)+(0,3)-form
$d\Omega$ on $M$  satisfies the identity \beq\label{nk3}
-\frac{1}{3}(\nabla^g_Xd\Omega)(Y,Z)=(\nabla^g_X\nabla^g_YJ)Z
=\frac{1}{2}a\left(g(Y,JZ)X+g(X,Z)JY-g(X,Y)JZ\right), \eeq where
$d\Omega(X,Y)$ denotes the (1,2) tensor corresponding to the
3-form $d\Omega$ via the metric $g$.

Since $d\Omega(X,Y)$ is normal to $L$ for $X,Y$ tangent to $L$, we
obtain $D_XJY=-\frac{1}{3}d\Omega(X,Y)+J\nabla^{g_L}_XY, \quad
A_{JY}X=-J\alpha(X,Y)$. Using these properties of the second
fundamental form, we calculate
\beq(\nabla^g_Xd\Omega)(Y,Z)=-A_{d\Omega(Y,Z)}X+D_Xd\Omega(Y,Z)-
d\Omega(\nabla^g_XY,Z)-d\Omega(Y,\nabla^g_XZ)= $$ $$
J\alpha(Jd\Omega(Y,Z),X)- \frac{1}{3}Jd\Omega(X,d\Omega(Y,Z))
-J(\nabla^{g_L}_XJd\Omega)(Y,Z)-d\Omega(\alpha(X,Y),Z)
-d\Omega(Y,\alpha(X,Z))
\label{stuffad} \eeq
 for
$X,Y,Z$ tangent to $L$.  Multiplying the last equality by $J$ and
using (\ref{nk3}) we get for the normal component \beq\label{las}
\alpha(Jd\Omega(Y,Z),X)+Jd\Omega(\alpha(X,Y),Z)
+Jd\Omega(Y,\alpha(X,Z))=0.
\eeq The last equality means $tr \alpha=0$. Indeed, we may assume
that $\frac{1}{3}Jd\Omega(e_1,e_2)=\sqrt{\frac{a}{2}}.e_3$  form
an orthonormal basis on $L$ for any even permutation of $(123)$.
Evaluating (\ref{las}) on those basis we get $tr \alpha=0$ by
taking the cyclic sum and using the skew-symmetry of $N$. Hence,
$L$ is a minimal submanifold of $M$. \hfill {\bf Q.E.D.}

 The theorem above generalizes the result of
Ejiri \cite{Ej} which states that for the Nearly K\"ahler $S^6$
any Lagrangian submanifold $M^3\subset S^6$ is minimal. SAS
calibrations in the Nearly K\"ahler $S^6$ will be consider below.

Since the torsion is a  $(3,0)+(0,3)$ form, the differential
system for the deformation of SAS calibrations  on a nearly
K{\"a}hler six-dimensional manifold $(M,J,g)$ reduces to the
equations
\beq\label{nksas}
dU=\frac{3}{4}(i_U)N=3(i_UT)=-(i_{JU})d\Omega, \quad \delta U=0\ .
\eeq

\begin{proposition}\label{tgm}

Let $V=JU$ be a SAS deformation of a SAS calibrated compact
submanifold $L$ on 6-dimensional Nearly K\"ahler manifold $M$.
Then the following formula holds
\beq\label{insas}
\int_L
(2Ric^g(U,U) - \frac{9}{2}a||U||^2 +\frac{1}{2}({\cal
L}_Ug)^2)dVol.|_L=0 \ . \eeq
In particular, if $M=S^6$ then $U$ cannot be
a Killing vector field on $L$.
\end{proposition}

{\it Proof:} We shall use the following general formula on a
compact Riemannian manifold  \cite{Yano}
\beq\label{ya}
\int_L(Ric_{ij}U^iU^j+(\nabla^jU^i)(\nabla_iU_j)-(\delta
U)^2)dVol.=0 \ .\eeq
The formula (\ref{ya}) follows from the identity
\beq
\nabla_i(U^j(\nabla_jU^i)-(\nabla_jU^j)U^i)=
Ric_{ij}U^iU^j+(\nabla^jU^i)(\nabla_iU_j)-(\nabla_jU^j)^2
\label{stuffae}
\eeq
by an integration over the compact $L$.

Let $V=JU$ be a SAS deformation. The constant type condition
(\ref{nk2}) implies $||dU||^2=9||i_UT||^2=9a||U||^2$. Substituting
the latter equality into (\ref{ya}) we get (\ref{insas}).

For any minimal lagrangian submanifold $L$ of the Nearly K\"ahler
$S^6$ we have $Ric^g(U,U)=a||U||^2 -
\sum_{i=1}^3g(\alpha(U,e_i),\alpha(U,e_i))$, where $e_1,e_2,e_3$
is an orthonormal basis on $L$. Substituting the last equality
into (\ref{insas}) and taking into account (\ref{nksas}), we get a
contradiction with the assumption that $U$ is Killing. \hfill {\bf
Q.E.D.}

\subsection{Hermitian manifolds with holonomy $SU(n)$}

In this section, we take $(M,J,g)$ to be a hermitian manifold,
 ${\rm dim} \ M=2n$, for
which the holonomy of either the Bismut connection $\nabla^b$
or the Chern connection $\nabla^c$ is contained in $SU(n)$.
The definitions of these connections are given in section fifteen.
Both these cases will emerge in the investigation of
hermitian manifolds
with trivial canonical bundle in sections sixteen and seventeen.

First consider the case for which the Bismut connection has holonomy
contained in $SU(n)$. In such a case the differential
system for the deformation of SAS
calibrations becomes
\bea
dU- {1\over2} U^a H_{a b' c'} J^{b'}{}_b J^{c'}{}_c e^b\land e^c&=&0
\nonumber \\
(\nabla^g)^aU_a-U^a \theta_a&=&0
\eea
where $H$ is the torsion of the Bismut connection and $\theta$ is
the Lee form (see section fifteen).
To derive the first equation, we have used the fact that the
torsion three-form $H$ of the Bismut
connection is (2,1) and (1,2) with respect to $J$; this
follows from the integrability
of the complex structure and the fact that $J$ is
parallel with respect to the Bismut connection.

There are two cases to consider. If the hermitian manifold
$(M,g,J)$ is balanced, then $\theta=0$, and the deformation
equations are
\bea
dU- {1\over2} U^a H_{a b' c'}
J^{b'}{}_b J^{c'}{}_c e^b\land e^c&=&0
\nonumber \\
(\nabla^g)^aU_a&=&0\ .
\eea
In particular, $U$ is co-closed. Next
assume that $(M,g,J)$ is conformally balanced, ie $\theta=2d f$
for some function $f$ on $M$. This class of hermitian structures
appears in some applications in physics. Rescaling $U=e^{2f} \hat
U$, we find that the differential system becomes
\bea
d\hat U+ 2 d
f \land \hat U - {1\over2} \hat U^a H_{ab'c'} J^{b'}{}_b
J^{c'}{}_c e^b\land e^c&=&0
\nonumber \\
(\nabla^{g})^a\hat U_a&=&0\ .
\eea
So $\hat U$ is again co-closed.

Next take $(M,J,g)$ to be a hermitian manifold for which the
associated Chern connection $\nabla^c$ has holonomy contained in
$SU(n)$. In this case the parallel (n,0)-form $\psi$ is
holomorphic and therefore closed, $d\psi=0$. The deformation
differential system (\ref{rewrb}, \ref{rewra}) for SAS calibrations
becomes
 \bea dU - U_a
(C^a+
  \hat C^a )&=&0
 \nonumber \\
(\nabla^g)^a U_a&=&0 \ ,
 \label{rewra1}
 \eea
where $C$ is the torsion of the Chern connection. Again,
 the one-form $U$ is co-closed.

\section{Examples}

\subsection{Group manifold examples}
\vskip 0.4cm
\leftline{{\underline {SAS calibrations in hermitian group manifolds}}}
\vskip 0.3cm

Consider the group manifold $S^3\times S^3$ with metric
\beq
ds^2=(\sigma^1)^2+(\sigma^2)^2+(\sigma^3)^2+(\tilde\sigma^1)^2
+(\tilde\sigma^2)^2+(\tilde\sigma^3)^2
\label{groupmet}
\eeq
and a complex
structure $J$ with associated K{\"a}hler form
\beq
\Omega=\sigma^1\land \sigma^2-\tilde \sigma^1\land \tilde \sigma^2
+\sigma^3\land\tilde\sigma^3
\label{groupforma}
\eeq
where $\sigma^i$ and $\tilde\sigma^i$ are left invariant one-forms
satisfying
\bea
d \s^a &=& -{1 \over 2} \e^{abc} \s^b \land \s^c \nonumber
\\
d \tilde\sigma^a &=& - {1 \over 2} \e^{abc} \tilde\s^b
 \land \tilde\s^c \ .
\eea
The associated parallel (3,0)-form is
\beq
\omega=e^{i{\pi\over4}}(\sigma^1+i\sigma^2)\land
 (\tilde\sigma^1-i\tilde\sigma^2)\land (\sigma^3+i\tilde\sigma^3)\ .
\eeq
The $X=S^3$ submanifold of $S^3\times S^3$ which is
 defined by the diagonal embedding
$\sigma^i|_X=\tilde \sigma^i|_X$ is a SAS calibration.

 This SAS calibration has moduli. To see this observe that
 both $\Omega$ and $\omega$ are
 invariant under the left action of $S^3\times S^3$. In
 addition $S^3\times S^3$ acts on the diagonal
 $S^3$ as $(k_1, k_2) (g,g)\to (k_1 g, k_2 g)$. Thus if the
 diagonal $S^3$ is a SAS calibration,
 then all the right cosets of $S^3$ in $S^3\times S^3$ are
  SAS calibrations as well. The moduli
 space of these deformations is $S^3\times S^3/S^3$. Observe
 that if $(k_1, k_2)=(k'_1 h, k'_2 h)$ for
 $h\in S^3$, then $(k_1, k_2)$  and $(k'_1, k'_2)$ generate the
 same deformation.

 For an alternative way to see this, let $\nabla$ be the
  connection on the group
 $S^3\times S^3$ associated with the left action.
 Observe that the metric and the K{\"a}hler
 form are parallel with respect to the connection $\nabla$.
 The connection $\tilde\nabla$ which has torsion $\tilde T=-T$
 is associated with the right action on the
 group manifold $S^3\times S^3$.
 In particular all right-invariant vector fields are parallel with
 respect to $\tilde \nabla$. Thus they satisfy
the equation (\ref{sasc}).
 This is equivalent to the analysis above in which the left
 group action was used. This is because the right-invariant vector fields
 generate the left-action on group manifolds.
Of course  the right-invariant vector fields which are tangent to
the diagonal $S^3$ generate diffeomorphisms of the diagonal $S^3$
and so they are not
 tangent to the moduli space. However, there are three
linearly independent right invariant
vector fields which are normal to the diagonal $S^3$,
which are given by
\beq
V^{(R)}_{(i)} = \rho_i - {\tilde{\rho}}_i
\label{righty}
\eeq
for $i=1,2,3$ where $\{\rho_i, \tilde \rho_i; i=1,2,3\}$
are right invariant
vector fields on  $S^3\times S^3$.
Hence the dimension
of the moduli space is at least three.

To find whether the moduli space has dimension
 more than three, one should
find the number of solutions to the differential
equations (\ref{sasa} and (\ref{sasb}) or equivalently (\ref{sasc}).
Adapting them to this example, we have
\bea
dU-U \land \s^3 + U_3 \s^1 \land \s^2 &=0 \nonumber
\\
\delta U - U_3 &=0
\label{ssmpeqn}
\eea
where $U \equiv U_i \s^i |_X$. After some
 computation, it can be shown
that the only solutions to ({\ref{ssmpeqn})} are
 given by linear combinations of
$i_{V^{(R)}_{(i)}} \Omega |_X$ as expected.
Hence the moduli space is three-dimensional.

We remark that there is a hermitian calibration which is a torus
$T^2=S^1\times S^1$ along the directions $(\sigma^3, \tilde
\sigma^3)$ of $S^3\times S^3$. This torus is the fiber of the
product of fibers of the product Hopf fibration $T^2\rightarrow
S^3\times S^3\rightarrow S^2\times S^2$. The homology class
$[T^2]$ is trivial because $H_2(S^3\times S^3)=0$. This is an
example of a family of hermitian calibrations with base space
$S^2\times S^2$.

A SAS calibration on the complex Iwasawa manifold is given at the
end of the paper, in the last section.

\vskip 0.4cm
\leftline{{\underline {SAS calibrations in almost
 hermitian group manifolds}}}
\vskip 0.3cm

{}For another group manifold example consider again the group
manifold $S^3\times S^3$ with metric (\ref{groupmet}) but now
equipped with the {\it almost} complex structure $J$ with
associated K{\"a}hler form
\beq \Omega=\s^1\land \tilde
\s^1+\s^2\land \tilde \s^2+\s^3\land \tilde \s^3\ . \eeq The
associated (3,0)-form is \beq \omega=(\s^1+i\tilde\s^1)\land
(\s^2+i\tilde\s^2)\land (\s^3+i\tilde\s^3)\ . \eeq
It is clear
that  the three-sphere given by $S^3\times \{e\}$ is a SAS cycle,
where $e$ is the identity element. It is also clear that any
three-sphere in $S^3\times S^3$ given by the embedding $S^3\times
\{k\}$, $k\in S^3$, is again a SAS cycle for the above generalized
calibration. Thus there is a moduli space which has dimension at
least three. This can also be derived using the connection
$\tilde\nabla$ as in the other example above. In fact the
dimension of the moduli space is three. To see this observe that
the deformation equations for SAS calibrations imply that
\bea
 \nabla_a U_b-\nabla_b U_a&=&0\nonumber\\
 \nabla_a U^a&=&0\ ,
 \eea
where $\nabla$ is the flat connection on $S^3$ with associated
frame the left-invariant 1-forms $\{\s^a; a=1,2,3\}$. These
equations in particular imply that
 \beq
 \nabla^a \nabla_a U_b=0\ .
 \label{stuffah}
 \eeq
 Using
 \beq
 \int_{S^3} ||\nabla U||^2=-\int_{S^3} (U,\nabla^2 U)=0
 \label{stuffai}
 \eeq
 we conclude that $U$ is parallel with respect to
 $\nabla$ and so left-invariant.
 Since there are three linearly independent
left-invariant vector fields on $S^3$,
  the dimension
 of the moduli space is three. In fact the moduli space is $S^3$.

We can also consider a similar group manifold
example as above but this time with
\beq
\omega=i(\s^1+i\tilde\s^1)\land (\s^2+i\tilde\s^2)
\land (\s^2+i\tilde\s^2)\
\eeq
as a (3,0) form. In this case a SAS cycle is $\{e\}\times S^3$.
 The moduli space is again
$S^3$.

We remark that in both the above group manifold examples,
there is an almost
hermitian calibration which is a torus $T^2$
along the directions $(\sigma^3, \tilde \sigma^3)$.
In fact the induced almost complex structure on $T^2$
 is integrable and so $T^2$ is complex.

The above two group manifold examples can be easily
 generalized as follows. Let
$G$ be a semisimple Lie group (${\rm dim}G=k$).
On the group manifold
$G\times G$, we can define the metric
\beq
ds^2(G\times G)=\sum_{a=1}^k \big((\s^a)^2 + (\tilde \s^a)^2\big)
\label{stuffaj}
\eeq
and the almost complex structure
$J$ with associated K{\"a}hler form
\beq
 \Omega=\sum_{a=1}^k \s^a\land \tilde\s^a
\label{stuffak}
 \eeq
where $\{\s^a; a=1,\dots, k\}$ and $\{\tilde\s^a; a=1,\dots, k\}$
are the left invariant one-forms of $G\times G$; the first set is
that of the first group in the product $G\times G$ while the
second set is that of the second group. We can also define a
(n,0)-form as
\beq \omega=(\s^1+i\tilde\s^1)\land\dots\land
(\s^k+i\tilde\s^k)\ .
\label{stuffal}
\eeq
The submanifold $G\times\{e\}$ is a SAS
calibration with respect to ${\rm Re}\ \omega$. In fact all spaces
$G\times \{h\}$, $h\in G$, are SAS calibrations. Therefore the
dimension of the moduli space is at least $k$. In fact it can be
shown that the dimension of the moduli
space  is exactly $k$ by repeating the analysis for
$S^3\times S^3$ examples above. In particular,
 it is straightforward to show
that the solutions of the deformation equations
 $U=-i_V \Omega |_X$ are left-invariant one-forms.
Similarly $\{h\}\times G$ are also
SAS calibrations with respect to ${\rm Re} \ \omega$ where
\beq
\omega=(i)^{-k} (\s^1+i\tilde\s^1)\land\dots\land
(\s^k+i\tilde\s^k)
\eeq
in this case. The dimension of the moduli
space is again $k$.

\subsection{SAS submanifolds in $S^3\times S^3$
 with the left-invariant Einstein metric}

There is another Nearly K\"ahler structure on $S^3\times S^3$ which
can be constructed as follows.
First write $Spin(4)=S^3\times S^3$ and decompose the
Lie algebra of $Spin(4)$ as $spin(4)=so(4)=h+m$ with
$h={\rm span}\{E_{12},E_{13},E_{23}\}$ and
$m={\rm span}\{E_{14},E_{24},E_{34}\}$, where
 the matrices $E_{ij},
(i<j)$ are the standard generators of $so(4)$.
 Denoting the associated left-invariant forms
as the elements of the basis, we have
\beq
dE_{ij}=-\sum_{k=1}^4E_{ik}\wedge E_{kj}\ .
\label{stuffba}
\eeq
Denote the Killing
form on $spin(4)$ by $B(X,Y)=-1/2 tr(XY)$. Then there are two
Einstein metrics on $S^3\times S^3$. One is associated with
the bi-invariant metric  $B_1=B|_{h\times
h}+B|_{m \times m}$. The SAS calibrations for
this manifold have already
been investigated above. The other is associated
 with the left-invariant metric
 $B_{1/3}=\frac{1}{3}B|_{h\times
h}+B|_{m \times m}$. With respect to $B_{1/3}$ we consider the
orthonormal basis $e_1=\sqrt{3}E_{12}, e_2=\sqrt{3}E_{13},
e_3=\sqrt{3}E_{23}, e_4=E_{14}, e_5=E_{24}, e_6=E_{34}$. In this
basis the K\"ahler form is
\beq
\Omega= -\frac{1}{2}\left(e^1\wedge
e^6-e^2\wedge e^5+ e^3\wedge e^4\right) .
\label{stuffbk}
\eeq
Denote the associated
almost complex structure with $J$. Then $(S^3\times
S^3,B_{1/3},J)$ is a Nearly K\"ahler non-K\"ahler manifold
\cite{Gru}.

Consider a copy of $S^3\subset S^3\times S^3$ determined by the
integrable distribution $h={\rm span}\{e_1,e_2,e_3\}$. This is a
lagrangian submanifold of the Nearly K\"ahler manifold $(S^3\times
S^3,B_{1/3},J)$ and therefore it is a SAS calibration.

The moduli space is at least 3-dimensional. Indeed, simple
calculations show that \beq de^j=-(i_{Je_j}d\Omega), \quad \delta
e^j=0 \label{stuffbc} \eeq
 for $ j=1,2,3$, where $e^j$ is the dual 1-form to
$e_j$. Hence, $e_1,e_2,e_3$ are solutions of the differential
system (\ref{nksas}).

\subsection{SAS calibrations on Flag manifold}
Let $F_{1,2}={\bf U}(3)/{\bf U}(1)\times {\bf U}(1)\times {\bf
U}(1)$ be the complex three-dimensional flag manifold. Consider
the reductive decomposition $ {\bf u}(3)={\bf h}\oplus{\bf m}$
where ${\bf u}(3)$ is the Lie algebra of the unitary group ${\bf
U}(3)$  and   $\bf h$ and $\bf m$ are determined by: ${\bf h}=
\cong {\bf u}(1)\oplus {\bf u}(1)\oplus {\bf u}(1)\subset {\bf
u}(3)$ and \beq
 {\bf m}= \left\{
\begin{array}{lll}
0       & a      & b       \\ -\bar{a}& 0      & c       \\
-\bar{b}&-\bar{c}& 0       \\
\end{array}
\right\} \subset {\bf u}(3). \label{stuffbe} \eeq Identifying  any
element $X \in TF_{1,2} \cong {\bf m}$ with the corresponding
triple of complex numbers $(a,b,c)$, we consider the {\bf
U(3)}-invariant Riemannian metric on $F_{1,2}, g(X,X)=\vert a\vert
^2+\vert b\vert ^2+\vert c\vert ^2$.
 An invariant almost complex structure on $F_{1,2}$ is defined by $ J:
(a,b,c) \rightarrow ( ia, -ib, ic) $ and it is compatible with the
invariant metric $g$. Then $(F_{1,2},g,J)$ is a Nearly K\"ahler
non K\"ahler 6-dimensional manifold. We consider an orthonormal
basis of $TF_{1,2}$ given by \beq
e_1=\frac{1}{\sqrt{2}}(1,0,0),\quad e_2=Je_1,\quad
e_3=\frac{1}{\sqrt{2}}(0,1,0),\quad e_4=-Je_3,\quad
e_5=\frac{1}{\sqrt{2}}(0,0,1),\quad e_6=Je_5.\label{stuffcb} \eeq
Then the 3-sphere $S^3$ determined by the integrable lagrangian
distribution $e_1,e_3,e_5$ is a SAS calibration on the Nearly
K\"ahler 6-dimensional $(F_{1,2},g,J)$.

\subsection{SAS calibrations and the 6-sphere}

Let ${\it Im}{\cal O}$ be the 7-dimensional vector space of
imaginary octonions. Consider the unit sphere $S^6 \subset {\it
Im}{\cal O}$. The right multiplication by  $u\in S^6$ induces a
linear transformation $J_u:{\cal O} \to {\cal O}$ which is
orthogonal and satisfies $J^2=-1$. The operator $J_u$ preserves
the 2-plane spanned by $1$ and $u$ and therefore preserves its
orthogonal 6-plane which may be identified with $T_uS^6$. Thus
$J_u$ induces an almost complex structure on $T_uS^6$ which is
compatible with the inner product induced by the inner product of
${\cal O}$. Therefore   $S^6$ has an almost complex structure
which is compatible with the standard metric $g$ on $S^6$ and so
$(S^6,g,J)$ is an almost hermitian manifold. In fact $(S^6,g,J)$
is Nearly K{\"a}hler.  The group of automorphisms is the
exceptional group $G_2$.

Let $L\subset S^6$ be a three-dimensional  Lagrangian submanifold
of $S^6$ with respect to a K{\"a}hler 2-form. Then by the results
of section 7.2, $L$ is a SAS calibration and is minimal. In
addition  it satisfies $g(\nabla^g_XJ)Y,Z)=0$ for $X,Y,Z$ tangent
vectors to $L$ by the result of N.Ejiri \cite{Ej}. For example
consider the invariant $G_2$ form in (\ref{gstring}) or
equivalently  in  (\ref{assocfrm}) below and view $\{e^i;
i=1,\dots,7\}$ an orthonormal basis in $\bR^7$. Then the K\"ahler
form at the point $x=x^ie_i$ of $S^6$, $\sum_{i=1}^7 (x^i)^2=1$,
is
\beq
\Omega= x^i \phi_{ijk} e^j\wedge e^k
\label{stuffcc}
\eeq
restricted in the
directions orthogonal to $x$. The three-sphere defined by the
equations $x^1=x^2=x^3=0$ is a Lagrangian submanifold and so a SAS
calibration.

\section{Deformations of generalized co-associative calibrations}

Let $(M,g, \psi)$ be a seven-dimensional manifold which admits a
metric connection $\nabla$ whose holonomy is contained in $G_2$.
As we have mentioned there is a local orthonormal frame $\{e^a \ ,
\ e^i \}$ for $a,b=4,5,6,7$ and $i,j=1,2,3$ such that the parallel
three-form $\psi$ takes the canonical form
\begin{equation}
 \psi = e^{123} + e^1 \land (e^{45}-e^{67}) +e^2 \land
(e^{46}+e^{57}) +e^3 \land (e^{47}-e^{56})\ , \label{assocfrm}
\end{equation}
where $e^{12}=e^1 \land e^2$ and similarly for the rest.
Observe that $\psi$ can also be written as
\beq
\psi = e^{123} +\sum_{i=1}^3 e^i \land \Omega_i
\label{assocfrma}
\eeq
where $\{\Omega_i, i=1,2,3\}$ is a basis of anti-self-dual two-forms
in the directions spanned by the $\{e^a; a=4,5,6,7\}$ frame basis.
The generalized co-associative calibrating four-form is simply the
Hodge dual of $\psi$, $* \psi$.

\begin{proposition}

A necessary and sufficient condition for a four-dimensional
submanifold $X$ of $M$ to be a generalized co-associative calibration
with respect to $* \psi$ is that $\psi|_X=0$.

\end{proposition}

{\it Proof:}
The proof of this proposition
is similar to that given for standard co-associative calibrations
in  \cite{harveylawson} and so it will not be repeated here.

\rightline{{\bf Q.E.D.}}

If $X$ is a co-associative calibrated submanifold,
 we can adapt a frame at
every point of $X$ such that the directions $\{e^a; a=4,5,6,7\}$
 are tangent to $X$ and $\{e^i; i=1,2,3\}$ are normal.
Expressing the condition ${\cal L}_V\psi=0$ for the deformation of
a co-associative calibration $X$ along the normal vector field $V$
in terms of the  $\tilde\nabla$ connection of (\ref{lie}), we have
\beq
(\Omega_i)_{ab} \tilde\nabla_c V^i e^a\land e^b\land e^c=0
\label{coco} \eeq
where $\{\Omega_i; i=1,2,3\}$ is the basis of
anti-self-dual K{\"a}hler forms used to construct the form $\psi$
in (\ref{assocfrm}). Using the anti-self duality of $\{\Omega_i;
i=1,2,3\}$, the equation (\ref{coco}) also implies
\beq
\sum_{i,a}(J_i)^a{}_b \tilde \nabla_a V^i=0\ , \label{cocona}
\eeq
where
$J_i$ are the (almost) complex structures associated with
$\Omega_i$.

There is an alternative way to express the deformation equations.
For this observe that the normal bundle $N_X$ of $X$ in $M$ is
isomorphic to the bundle $\Lambda^{2-}(X)$ of anti-self-dual
two-forms of $X$, $N_X=\Lambda^{2-}(X)$. The proof of this is
similar to that given by  \cite{ML}. It is based on the
observation that the normal bundle $N_X$ and  $\Lambda^{2-}(X)$
are both associated to the principal $SO(4)$ frame bundle of $X$
with the same representation, i.e. the three-dimensional
anti-self-dual representation of $SO(4)$. Note that the
seven-dimensional $\rho_7$ representation of $G_2$ which leaves
three form $\psi$  invariant decomposes as $\rho_7=v_3\oplus v_4$
under the action of $SO(4)\subset G_2$, where $v_3$ is the
three-dimensional anti-self-dual representation of $SO(4)$ acting
on the directions $123$ and $v_4$ is the standard four-dimensional
vector representation of $SO(4)$ acting on the directions $4567$.
Let $V=V^i e_i$ be a normal vector field of $X$,  then the isomorphism
is given by $\a_V=i_V \psi |_X$. Observe that $\a_V$ is an
anti-self-dual two-form on $X$.

Next the condition ${\cal L}_V\psi=0$ can be written as
\begin{equation}
d \a_V + i_V d \psi |_X=0\ . \label{stuffbb}
\end{equation}
If   $d \psi =0$, then the dimension of the moduli space
is equal to $b_-^2(X)$, i.e. the dimension of the space of
 anti-self-dual
harmonic two-forms of $X$.
Now we shall turn to the case where $d \psi \not=0$.
Using the fact
that $\psi$ is parallel with respect to
$\nabla$, we find
that
\beq
d (\a_V)- (\a_V)_{ab} T^a \land e^b +{1 \over 2} \psi_{ibc}
V^j T^i{}_{ja} e^a \land e^b \land e^c=0 \ .
\eeq
Using the relation
\begin{equation}
g_{AB} = {1 \over 6} \psi_{A CD} \psi_B{}^{CD}\ , \label{giimet}
\end{equation}
 the deformation equations may be written solely in terms of $\a_V$
as
\beq
d \a_V + (\a_V)_{ab} \big[ - T^a \land e^b
+{1 \over 6} \psi_{icd} \psi^{jab}
T^i{}_{j a_1} e^c \land e^d \land e^{a_1} \big] =0\ .
\label{coassocsd}
\eeq

Furthermore, as $\a_V$ is anti-self-dual, this defines an elliptic
system of partial differential equations. Hence we conclude that

\begin{corollary}
The moduli space of generalized co-associative calibrations,
 if it exists,  is finite dimensional.
\end{corollary}

\section{Deformations of generalized associative calibrations}

Let $(M,g, \psi)$ be a seven-dimensional manifold which admits a metric
connection $\nabla$ whose holonomy is contained in $G_2$. Such
manifolds also admit  generalized
associative calibrations, in addition to
the generalized co-associative calibrations
investigated in the previous section. The former is a degree
 three calibration
associated with the  three form $\psi$ of (\ref{assocfrm}). We
again introduce the orthonormal frame $\{e^A\}=\{e^i \ , \ e^a;
i=1,2,3, a=4,5,6,7 \}$, where now  $\{e^i \}$ span the tangent
directions of the generalized associative cycle and $ \{ e^a \}$
span the normal directions.

As in the case of standard calibrations,
the condition for a three-dimensional submanifold $X$ to be
calibrated with respect to $\psi$ is that a certain
vector-valued three-form $\chi\in \Omega^3(M, TM)$ should
vanish on $X$. The
form $\chi$ is related to the cross product on $\bR^7={\rm Im} \bO$
 and it is
invariant under $G_2$, so $\chi$ is $\nabla$-parallel. In particular
in the basis that we have written
the three-form $\psi$, $\chi$ is given by
\beq
\chi = \sum_{A=1}^7\chi^A \otimes e_A
\label{stuffce}
\eeq
where
\begin{eqnarray}
\chi^1 &=& \big(e^{256}-e^{247}+e^{346}+e^{357}\big)  \cr \chi^2 &=&
\big(e^{147}-e^{156}-e^{345}+e^{367}\big)  \cr \chi^3 &=&
\big(e^{245}-e^{267}-e^{146}-e^{157}\big)  \cr \chi^4 &=&
\big(e^{567}-e^{127}+e^{136}-e^{235}\big)  \cr \chi^5 &=&
\big(e^{126}-e^{467}+e^{137}+e^{234}\big)  \cr \chi^6 &=&
\big(e^{457}-e^{125}-e^{134}+e^{237}\big)  \cr \chi^7 &=&
\big(e^{124}-e^{456}-e^{135}-e^{236}\big)\ . \label{vannish}
\end{eqnarray}

To compute the dimension of the moduli space of generalized associative
calibrations, we require that $\cL_V \chi |_X =0$. Using the fact
that this vector-valued three-form is parallel with respect to the
connection $\nabla$, (\ref{lie}) and after observing that $(\chi^A
\cL_V e_A)|_X=0$, we find that $\cL_V \chi |_X =0$ implies
\beq
\sum_{i,b} (\Omega_i)_{ab} \tilde\nabla_i V^b=0 \label{diraca}
\eeq
where $\tilde\nabla$ is the connection with torsion $\tilde
T=-T$. The normal bundle $N_X$ of a generalized associative
submanifold is isomorphic to the spin bundle $\bS$ of $X$. This
can be shown by observing that both $N_X$ and $\bS$ are associated
to the Spin principal bundle, $\tilde P$, of $X$, which is the
double cover of the frame bundle of $X$, with the same
representation. The proof of this is similar to that for standard
 calibrations and it has been described in \cite{harveylawson, ML}. Here
 we shall summarize the proof.  First observe that  every oriented
 three manifold admits a spin structure and so $\tilde P$
 exists for all associative calibrations $X$.
Then observe that the seven-dimensional representation
 $\rho_7$ of $G_2$ which leaves
  the three-form $\psi$ invariant decomposes as
 $\rho_7=v_3\oplus s_4$, where
 $v_3$ is the three-dimensional representation
 of $SU(2)$ induced by the
 standard three-dimensional vector representation
of $SO(3)$ acting on
the directions $123$ and $s_4$ is the four-dimensional real spinor
representation
 of $SU(2)$ acting on the directions $4567$.
 Since $4567$ are the normal
 directions of $X$, the normal bundle $N_X$
and the spin bundle $\bS$
 are associated to $\tilde P$ with the same
 representation $s_4$, so $N_X=\bS$.
Therefore the deformation equation (\ref{diraca}) is the Dirac
equation in three-dimensions with respect to the connection
$\tilde\nabla$;
 the gamma-matrices
are given by the $\{\Omega_i; i=1,2,3\}$. This
 is in fact an elliptic
differential equation and so if a moduli exist,
 the moduli space is finite dimensional.

The index of the Dirac operator that appears in the deformations
of generalized associative calibrations vanishes. Because of this
it is expected that generic generalized associative calibrations
will not have moduli. This is similar to the case of standard
associative calibrations. Although generic generalized associative
calibrations do not have moduli, we shall find many examples of
families of generalized associative calibrations in special cases.

\section{Special Cases}

There are several special cases of $G_2$ structures according to
Fernandez-Gray classification \cite{FG} depending on various
additional conditions that the three-form $\psi$ and its dual
$* \psi$ satisfy.
 \vskip 0.3cm
\leftline{\underline{Calibrated and cocalibrated $G_2$ manifolds}}
\vskip 0.3cm
The manifold $(M,g,\psi)$ is calibrated if $d\psi=0$
and $(M,g,\psi)$ is cocalibrated if $d*\psi=0$. It is known that
if $(M,g,\psi)$ is both calibrated and cocalibrated, then the
holonomy of the Levi-Civita connection $\nabla^g$ is contained in
$G_2$.

For calibrated $G_2$ manifolds, the generalized associative
submanifolds are minimal because $\psi$ is closed. In addition,
the deformation equations of such submanifolds are given by the
Dirac equation in (\ref{diraca}) with respect to a connection of
the normal bundle which is induced from a connection on $M$ which
has non-vanishing torsion.

For calibrated $G_2$ manifolds, the generalized co-associative
submanifolds generically are not minimal because $*\psi$ is not
closed. The deformations of such submanifolds are given in
(\ref{coco}) or in (\ref{stuffbb}). Because $d\psi=0$, the latter
equation can be simplified to
\beq
d\alpha_V=0\ .
\label{stuffda}
\eeq
Since
$\alpha_V$ is anti-self-dual, the dimension of the moduli space of
generalized co-associative calibrations in calibrated $G_2$
manifolds is $b^-_2$.

For cocalibrated $G_2$ manifolds, the generalized associative
submanifolds are not generically minimal because $\psi$ is not
closed. The deformation equations of such submanifolds  are given
by the Dirac equation (\ref{diraca}) with respect to a connection
of the normal bundle which is induced from a connection on $M$,
which has non-vanishing torsion.

For cocalibrated $G_2$ manifolds, the generalized co-associative
submanifolds  are  minimal because $d*\psi=0$. The deformations of
such submanifolds are given in  (\ref{coco}) or in
(\ref{stuffbb}).

Another type of $G_2$ manifold for which the associated
generalized calibrations can be analyzed as for co-calibrated
$G_2$ manifolds is that of cocalibrated $G_2$ manifolds of pure
type. For such manifolds $d*\psi=0$ and $d\psi\wedge \psi=0$.
Again the co-associative calibrations are minimal.

\vskip 0.3cm

\leftline{\underline{Integrable $G_2$ manifolds}}

\vskip 0.3cm

The manifold $(M,g, \psi)$ is an integrable  $G_2$ manifold iff
\beq
d*\psi=\theta \wedge *\psi
\label{stuffdb}
\eeq
where
$3\theta=-*(*d\psi\wedge \psi)$ is the Lee form.
It has  been shown in \cite{FI} that such $G_2$ manifolds admit a
unique connection with torsion a three-form.

For generic integrable $G_2$ manifolds, both $\psi$ and $*\psi$
are not closed, so the generalized associative and co-associative
calibrations are not minimal. The deformation equations for
generalized associative calibrations are given by the Dirac
equation (\ref{diraca}) but in this case the connection on the
normal bundle is induced from a connection with torsion a
three-form on $M$. The deformation equations for generalized
co-associative calibrations (\ref{stuffbb}) can be simplified
somewhat using the expression for the torsion in \cite{FI}. In
particular, denoting the torsion 3-form by $T=H$, we have
\cite{FI} that
\beq
d \psi = {1 \over 6} (d \psi . * \psi) * \psi + \theta
\land \psi + * H \ . \label{firstd}
 \eeq
Hence, noting that $i_V * \psi |_X =0$, it follows that
\beq
i_V d
\psi |_X = \big( - \theta \land \a_V +i_V (* H) \big) |_X
\label{ttska}
\eeq
and so we require that
\beq
d \a_V - \theta \land \a_V +{1 \over 36} (\a_V)_{a_1 a_2}
 \psi^{i a_1 a_2} *
H_{iabc} e^a \land  e^b \land  e^c =0\ . \label{ttskb}
\eeq
However, despite this simplification, it has not been possible to
compute the dimension of the moduli space.

The subclass of integrable $G_2$ manifolds which have applications
in physics (string theory) are those for which the 1-form $\theta$ is
exact and so $\theta = -2 d \Phi$ for some function $\Phi$ on $M$
which is identified with the dilaton. The analysis of generalized
associative and co-associative calibrations in this case is as for
the integrable $G_2$ manifolds above. There is some additional
simplification though in the deformation equations for generalized
co-associative calibrations. In particular, defining $\rho_V= e^{2
\Phi} \alpha_V$ and substituting in (\ref{ttskb}), we obtain
 \beq
d \rho_V +{1 \over 36} (\rho_V)_{a_1 a_2} \psi^{i a_1 a_2} *
H_{iabc} e^a \land  e^b \land  e^c =0\ . \label{ttskb1} \eeq

\vskip 0.3cm
\leftline{\underline{Nearly parallel or weak holonomy $G_2$ manifolds}}
\vskip 0.3cm

The manifold $(M,g,\psi)$ admits a nearly parallel or weak
holonomy $G_2$ structure iff $d\psi=\lambda *\psi$, for
 $\lambda$ constant.
If $\lambda=0$, then $(M,g,\psi)$ is calibrated. If
$\lambda\not=0$, then $(M,g,\psi)$ is co-calibrated. Since we have
already investigated the case of calibrated $G_2$ manifolds, we
shall focus on the case that $\lambda\not=0$. In \cite{FI}, it has
been shown that  nearly parallel manifolds admit a connection
$\nabla$ with torsion a three-form. In particular $  T=-{1\over
6}\lambda \psi$.

There are many examples of nearly parallel $G_2$ manifolds which
include $S^7$, $SO(5)/SO(3)$ and the Aloff-Wallach spaces
$N(n,m)=SU(3)/U(1)_{n,m}$;  the embedding
 of $U(1)$ in $SU(3)$ will be described later.

For nearly parallel $G_2$ manifolds, the generalized associative
submanifolds are not generically minimal because $\psi$ in not
closed. The deformation equations of such submanifolds  are given
by the Dirac equation (\ref{diraca}). Using the connection with
the torsion the three-form which is proportional to $\psi$, the
deformation equation can be simplified to
\beq
\sum_{i,b}
(J_i)^a{}_b \nabla_i^g V^b+{\lambda\over4} V^a=0\ .
\label{stuffdc}
\eeq
Therefore
the deformations of the associative submanifolds are eigenspinors
of the Dirac operator.

For nearly parallel $G_2$ manifolds, the generalized
co-associative  submanifolds  are  minimal because $d*\psi=0$. The
deformations of such submanifolds are given in  (\ref{coco}) or in
(\ref{stuffbb}). The latter equation can be simplified. In
particular using the definition of the nearly parallel $G_2$
manifold, we find that
\beq
d\alpha_V=0
\label{stuffdd}
\eeq
and therefore since $\alpha_V$ is anti-self-dual, it is harmonic.
The dimension of the moduli space of a generalized co-associative
calibration $X$ in a nearly parallel $G_2$ manifold $(M,g,\psi)$
is $b_2^-(X)$.

\section{Examples}

\subsection{A group manifold example}

Consider the group manifold $M=S^3\times \tilde S^3\times S^1$
with left-invariant metric
\beq
 g=\sum_i (\s^i)^2+\sum_i
(\tilde\s^i)^2+ (\s^0)^2
\label{stuffde}
\eeq
and equipped with the left-invariant
three-form
\beq
 \psi=\s^{123}+\s^1\wedge (\tilde \s^{01}-\tilde
\s^{23})+\s^2\wedge (\tilde \s^{02}+\tilde \s^{13})
+\s^3\wedge(\tilde \s^{03}-\tilde \s^{12})\ ,
\label{stuffdh}
\eeq
where $\{\s^i;
1,2,3\}$ and  $\{\tilde\s^i; 1,2,3\}$ are the left-invariant
one-forms on the three-spheres $S^3$ and  $\tilde S^3$ in $M$,
respectively and $\tilde\s^0$ is the invariant one-form on $S^1$.
Clearly this three-form $\psi$ defines a $G_2$ structure on $M$
which is parallel with respect to the $\nabla$-connection on the
group manifold associated with the left action.

It can be easily seen that the submanifold $S^3$ is a generalized
associative calibration, while  $S^1\times \tilde S^3$ is a
generalized co-associative calibration.

Observe that the submanifolds $S^3\times \{p\}$, $p\in \tilde
S^3\times S^1$ are all generalized associative calibrations and so
the moduli space has dimension at least four. In fact the moduli
space has dimension exactly four. To see this observe that the
equation for the deformations in this case is

\beq
 \sum_{i,b}(\Omega_i)_{ab} \tilde \nabla_i
V^b=\sum_{i,b}(\Omega_i)_{ab}  \nabla_i V^b=0
\label{stuffdi}
\eeq
 and $\nabla$ is
a flat connection. Therefore
\beq
 \nabla^2 V^a=0\ .
\label{stuffdj}
\eeq
Then
\beq
\int_{S^3} (\nabla V, \nabla V)=- \int_{S^3} (V, \nabla^2 V)=0
\label{srtuffdk}
\eeq
and hence $V$ is constant. So the moduli space has dimension four.
In fact the moduli space in this case is $\tilde S^3\times S^1$
and therefore $M$ is a family of generalized associative
calibrations.

Similarly, observe that the submanifolds $\tilde S^3\times
S^1\times \{p\}$, $p\in S^3$, are all generalized co-associative
calibrations and so the moduli space in this case has dimension at
least three. In fact the moduli space has dimension exactly three.
The deformation equation is
\beq
 \sum_{i,b}(J_i)^b{}_a \tilde
\nabla_b V^i=\sum_{i,b}(J_i)^b{}_a  \nabla_b V^i=
\sum_b \nabla^b(\alpha_V)_{ba}=0\ .
\label{stuffdl}
\eeq
Therefore $\alpha_V$ is co-closed
with respect to the flat connection $\nabla$. Since $\alpha_V$ is
anti-self-dual, it is also closed, $\nabla\wedge \alpha_V=0$,
with respect to $\nabla$. Since $\alpha_V$ is both closed and
co-closed is harmonic with respect to the Laplacian $\nabla^2$. A
partial integration argument similar to the one above implies that
$\alpha_V$  is necessarily  $\nabla$-parallel. This implies that
the dimension of the moduli space is three. In fact in this case,
the moduli space is $S^3$.

\subsection{Generalized associative calibrations in homogeneous spaces }

\vskip 0.3cm

\leftline{\underline { Generalized calibrations in $S^7=Sp(2)/Sp(1)$}}

\vskip 0.3cm

Identify $\bR^8=\bH^2$. Then observe that the action of $Sp(2)$
preserves the equation for $S^7$ written in terms of quaternions
with stability subgroup $Sp(1)$ up to a conjugation. This implies
that $Sp(2)/Sp(1)=S^7$. In addition observe that $Sp(1)\subset
Sp(1)\times Sp(1)\subset Sp(2)$. This leads to the principal
fibration $Sp(1)\rightarrow S^7\rightarrow \bH P^1$. This is the
principal fibration associated with the anti-self dual
$SU(2)=Sp(1)$ instanton connection in $S^4=\bH P^1$. Let
$\{\alpha^i; i=1,2,3\}$ be the associated connection with
curvature
\beq
 \omega^i=d\alpha^i+\epsilon^i{}_{jk} \alpha^j\wedge \alpha^k\ .
\label{stuffdm}
\eeq
The Bianchi identity implies that
\beq
d\omega^i=2\epsilon^i{}_{jk}\omega^j\wedge \alpha^k\ .
\label{stuffdn}
\eeq
In addition there is a local frame $\{ \ell^a; a=4,\dots,7\}$ such
that
\beq
\omega^i={1\over2} \Omega^i_{ab} \ell^a\wedge \ell^b\ ,
\label{stuffdo}
\eeq
where $\{\Omega^i; i=1,2,3\}$ is the basis of constant
anti-self-dual two-forms in $\bR^4$ given in (\ref{assocfrm}) and
(\ref{assocfrma}).

Next consider the metric and the three-form on $S^7$
\bea ds^2&=&
y^2 \sum_{i=1}^3 (\alpha^i)^2+ z^2 \sum_{a=4}^7 (\ell^a)^2 \cr
\psi&=& y^3 \alpha^1\wedge \alpha^2\wedge \alpha^3+ y z^2
\sum_{i=1}^3 \alpha^i\wedge \omega^i\ , \label{gtwosphere}
\eea
where $y,z\in \bR-\{0\}$. It can be easily seen by setting
$\{e^i=y \alpha^i; i=1,2,3\}$ and $\{e^i=z \ell^i; i=4,\dots,7\}$
that the metric $ds^2$ and $\psi$ above take the canonical form of
a $G_2$ structure as in (\ref{assocfrm}) and (\ref{assocfrma}).

The fibres of the fibration $Sp(1)\rightarrow S^7\rightarrow \bH
P^1$ are all associative generalized calibrations; this can easily be
 seen by observing that
\beq
\psi|_{Sp(1)}=y^3( \alpha^1\wedge
\alpha^2\wedge \alpha^3)|_{Sp(1)}= d{\rm vol}(Sp(1))\ .
\label{stuffdp}
\eeq
This is
the case for any $y,z\in \bR-\{0\}$. Therefore this fibration is a
family of generalized associative calibrations.

It can be easily seen that the $G_2$ structure on $S^7$ in
(\ref{gtwosphere}) is nearly parallel, $d\psi=\lambda * \psi$,
iff
\bea
-3y&=&\lambda z^2 \cr {1\over2}y^2+z^2&=&-{1\over2}
\lambda y z^2\ .
\eea
This system has a solution for
$y=-3/\lambda$ and $z=\pm 3/\lambda$. This gives a nearly parallel
$G_2$ manifold which is  the squashed $S^7$. Clearly the squashed
$S^7$ is a family of generalized associative calibrations.

\begin{remark}

The Hopf fibration $S^1\rightarrow S^3\rightarrow S^2$ is also a
smooth family of generalized calibrations. To see this observe
that the metric on $S^3$ can be written as
\beq
ds^2(S^3)= (\s^3)^2+ (\s^1)^2+(\s^2)^2
\label{stuffdq}
\eeq
  where $\s^1,
\s^2,\s^3$ are the left-invariant one-forms on $S^3$ and
$ds^2(S^2)=(\s^1)^2+(\s^2)^2$. It can be easily seen that $\s^3$
is a generalized calibration in $S^3$ of degree one. The
calibrated lines are circles which are the fibres of the
Hopf-fibration. Therefore $S^3$ is a family of generalized degree
one calibrations with space of parameters $S^2$.

\end{remark}

\vskip 0.3cm

\leftline{\underline {Generalized calibrations in $M=SO(5)/SO(3)$}}

\vskip 0.3cm

We shall demonstrate that $M=SO(5)/SO(3)$ is a family of
generalized associative calibrations. We remark that $M$ is not
homeomorphic to $S^7$; $M$ and $S^7$ have the same deRham
cohomology but $M$ exhibits torsion in the third cohomology.
Observe that $so(5)=so(4)\oplus \bR^4$ and $so(4)$ acts with the
fundamental representation on $\bR^4$. Since $so(4)=so(3)\oplus
so(3)$, the structure constants decompose under the decomposition
$\Lambda\bR^4=\Lambda^{2+}\bR^4\oplus \Lambda^{2-}\bR^4$. Under
this decomposition of $so(5)$ a frame can be introduced at $M$
which satisfies the following structure equations:
 \bea
de^a&= &(J_i)^a{}_b \rho^i\wedge e^b+ (I_i)^a{}_b
 \sigma^i\wedge e^b
\\
d\rho^i&=& \epsilon^i{}_{jk}\rho^i\wedge
\rho^j-{1\over2}(\Omega_i)_{ab} e^a\wedge e^b\ ,
 \eea
 where $\{e^a; a=1, \dots, 4\}$ are associated with a basis in
$\bR^4$, $\{\sigma^i; i=1,2,3\}$ are associated with a basis in
the Lie algebra of the stability subgroup of the coset and
$\{\rho^i; i=1,2,3\}$ are the rest of the generators. The
structure constant matrices $\{J_i; i=1,2,3\}$ are anti-self-dual
and the structure constants $\{I_i; i=1,2,3\}$ are self-dual;
$(\Omega_i)_{ab}=\delta_{ac} (J_i)^c{}_b$.

The metric and three-form on $M$ are the following:
\bea
ds^2(S^7)&=& \sum_{i=1}^3 (\rho^i)^2+\sum_{a=1}^4 (e^a)^2
\\
\psi&=&\rho^1\wedge \rho^2\wedge \rho^3+2
\sum_{i=1}^3(\Omega_i)_{ab} \rho^i\wedge  e^a\wedge e^b\ .
\label{datat}
\eea
Observe that both are invariant under the
stability group $SO(3)$ of the coset and so they are globally
defined on $M$. It can be easily seen that the data given in
(\ref{datat}) define a nearly parallel $G_2$ structure on $M$.

Next consider the obvious subgroups of $SO(5)$; $SO(3)\subset
SO(4)\subset SO(5)$. Then there is a fibration
$SO(4)/SO(3)\rightarrow SO(5)/SO(3)\rightarrow SO(5)/SO(4)$ or
equivalently $S^3\rightarrow M\rightarrow S^4$. The cotangent
bundles of the fibres at every point $p\in S^4$ are spanned by
$\{\rho^i|p; i=1,2,3\}$ and because the metric on $M$ is diagonal
in this basis all the $S^3$ fibres of this fibration are
associative submanifolds of $M$. This  fibration is a family of
generalized associative $G_2$ calibrations.

\vskip 0.3cm

\leftline{\underline {Generalized calibrations in Aloff-Wallach
 spaces $N(n,m)$}}

\vskip 0.3cm

Another class of nearly parallel $G_2$ manifolds are the so called
Aloff-Wallach spaces $N(n,m)= SU(3)/U(1)_{n,m}$. The $U(1)$ is
embedded in $SU(3)$ as
 \beq
{\rm diag} (e^{i n\chi}, e^{i m \chi}, e^{-i(n+m)\chi})\ ,
\label{uonee}
\eeq
 where $n,m\in\bZ$.
To construct the $G_2$ structures on this space write
$su(3)=u(1)\oplus \bR^7$. Under the action of $U(1)$, $\bR^7$
decomposes as $\bR^7=\bR^2\oplus \bR^2\oplus \bR^2\oplus \bR$.
This can be seen by using the action of $u(1)$ on the Cartan
subalgebra and the step operators of $su(3)$. In particular, each
$\bR^2$ is spanned by the step operators $E_{\pm \alpha}$, where
$\alpha$ is a positive root, while $\bR$ is spanned by the
direction in the Cartan subalgebra of $su(3)$ which is orthogonal
to the generator of the embedded $u(1)$. A local frame can be
introduced on $N(n,m)$ according to this decomposition as
$\{\sigma^i, \rho^i, \zeta^i, \eta; i=1,2\}$.  To be precise, let
$\{L_B{}^A; A,B=1,2,3\}$, $(L_B{}^A)^\dagger= L_B{}^A$, ${\rm tr}
L_A{}^A=0$, be the left invariant forms on $SU(3)$, $dL_A{}^B=i
L_a{}^C\wedge L_C{}^A$. We set $\sigma^1+i\sigma^2=i L_3{}^1$,
$\rho^1+i\rho^2=L_2{}^3$ $\zeta^1+i\zeta^2= L_1{}^2$, and
$\eta=\sqrt{2} (\cos\delta L_1{}^1+\sin\delta L_2{}^2)$, where
${\rm tan} \delta=-{n\over m}$. This decomposition is similar
 to that in
\cite{cvet}.  A metric on $N(n,m)$ can be
written as
 \beq ds^2=x^2 \sum_{i} (\sigma^i)^2+y^2 \sum_{i} (\rho^i)^2+z^2
\sum_{i} (\zeta^i)^2+ f^2 \eta^2\ , \label{metmetn} \eeq where
$x,y,x,f\in\bR-\{0\}$. To define a $G_2$ three-form, it is most
convenient to induce it from a K\"ahler form $\omega$ and a
(3,0)-form $\phi$ on $\bR^6$ because of the above decomposition of
$\bR^7$. Indeed consider the two-form
\beq
\omega= x^2
\sigma^1\wedge \sigma^2+y^2  \rho^1\wedge \rho^2+z^2
\zeta^1\wedge \zeta^2
\label{stuffha}
\eeq
and the (3,0)-form
\beq
 \phi= xyz
(\sigma^1+i\sigma^2)\wedge (\rho^1+i\rho^2)\wedge
(\zeta^1+i\zeta^2)\ .
\label{stuffhb}
\eeq
 Then the $G_2$ three-form on $N(n,m)$ can
be defined as \beq \psi={\rm Re}\phi- f \eta\wedge \omega\ .
\label{ffn} \eeq Setting $e^1=x\sigma^1, e^5=x\sigma^2,
e^2=y\rho^1, e^6=y\rho^2, e^3=z \zeta^1, e^7=z \zeta^2,
e^4=f\eta$, one can bring  the metric $ds^2$ and $\psi$ above into
 the canonical form of a $G_2$ structure given in
(\ref{assocfrm}) and (\ref{assocfrma}).

For all $x,y,z,f\in\bR-\{0\}$, the above data define a $G_2$
structure on $N(n,m)$. However not all these $G_2$ structures are
nearly parallel. It can be shown that if
 \bea
\lambda&=&(x^2+y^2+z^2)
\cr
4xyz+2 \sqrt{2} f \big(y^2 (\cos\delta-\sin\delta)+z^2
\sin\delta\big)&=& \lambda y^2 z^2
\cr
4xyz+2 \sqrt{2} f \big(x^2 (\cos\delta-\sin\delta)-
z^2 \cos\delta\big)&=& \lambda x^2 z^2
\cr
4xyz+2 \sqrt{2} f \big(x^2 \sin\delta-y^2
\cos\delta\big)&=& \lambda y^2 x^2
\eea
then the $G_2$ structure is nearly parallel.
It is known that these equations have solutions and so
there are nearly parallel $G_2$ structures on $N(n,m)$;
 for a recent discussion see \cite{swann, cvet}.

To find generalized $G_2$ calibrated submanifolds in $N(n,m)$,
 observe that
$U(1)_{n,m}\subset S(U(2)\times U(1))\subset SU(3)$.
Viewing $S(U(2)\times U(1))$ as a $3\times 3$
matrix, the embedding of $U(1)$ in $S(U(2)\times U(1))$ is as
 in (\ref{uonee}).
This sequence of subgroups of $SU(3)$ define the fibration
\beq
S(U(2)\times U(1))/U(1)_{n,m}\rightarrow N(n,m)
\rightarrow \bC P^2\ .
\label{stuffhc}
\eeq
 In fact it turns out that
the typical  fibre is $S(U(2)\times U(1))/U(1)_{n,m}=S^3/\bZ_p$,
where $p=|n+m|>0$; for $p=0$ the typical
fibre is $S^2\times S^1$.

Next decompose $s(u(2)\oplus u(1))$ under the action of $u(1)$ as
$s(u(2)\oplus u(1))=u(1)\oplus \bR^3$. Moreover $\bR^3$ decomposes
under the irreducible two dimensional real representation of
$u(1)$ as $\bR^3=\bR^2\oplus \bR$. Since $s(u(2)\oplus
u(1))\subset su(3)$, $\bR^3$ is a submodule of $\bR^7$ under the
action of $u(1)$. Therefore it can be arranged such that
 the tangent space of the fibres of the fibration
 $S^3/\bZ_p\rightarrow N(n,m)
\rightarrow \bC P^2$ is spanned by $\{\eta, \sigma^1, \sigma^2\}$.
It is clear that $\psi|_{S^3/\bZ_p}= d{\rm vol}(S^3/\bZ_p)$
and so every fibre
is a generalized associative calibration.
Therefore the Aloff-Wallach spaces
are families of generalized associative calibrations
for any $G_2$ structure
defined in (\ref{metmetn}) and (\ref{ffn}).

\section{Deformations of generalized Cayley Calibrations}

Let $(M,g)$ be an eight-dimensional Riemannian manifold which admits
a metric connection $\nabla$ with holonomy contained in $Spin(7)$.
On such a manifold there is a local frame
$\{e^A; A=1,\dots,8\}$ such that
the self-dual four-form
\begin{equation}
\Phi = e^{1234} + (e^{12}-e^{34}) \land (e^{56}-e^{78})
 +(e^{13}+e^{24}) \land (e^{57}+e^{68})
+(e^{14}-e^{23}) \land (e^{58}-e^{67}) +e^{5678}
\label{caylfrm}
\end{equation}
is $\nabla$-parallel.

The condition that a four-dimensional submanifold
 $X\subset M$ is calibrated with respect
 to $\Phi$ is that $\tau|_X=0$ where $\tau\in \Omega^4(M, F)$ is a
four-form which takes values on the vector bundle $F$;
 $F=P\times_{d_7} \bR^7$
where $d_7$ is the
seven-dimensional representation of $Spin(7)$, i.e.
 the one induced from the
standard seven-dimensional vector representation of
$SO(7)$. This four-form
$\tau$ is associated with the four-fold cross product
 of $\bR^8=\bO$ with values
in ${\rm Im}\bO$ and it is $Spin(7)$ invariant, so
 $\tau$ is $\nabla$-parallel
 Explicitly, in an appropriate basis, $\tau$ is
\beq
 \tau = \pmatrix{
(e^{14}-e^{23})\land (e^{57}+e^{68}) - (e^{13}+e^{24})\land (e^{58}
-e^{67})
\cr
(e^{12}-e^{34})\land
(e^{58}-e^{67}) - (e^{14}-e^{23})\land (e^{56}-e^{78})
\cr
(e^{13}+e^{24})\land (e^{56}-e^{78}) -
(e^{12}-e^{34})\land (e^{57}+e^{68})
 \cr
 e^{2345}-e^{1346}+e^{1247}-e^{1238}+e^{1678}-e^{2578}+e^{3568}-e^{4567}
\cr
e^{2346}+e^{1345}+e^{1248}+e^{1237}-e^{2678}-e^{1578}-e^{4568}-e^{3567}
\cr
e^{2347}+e^{1348}-e^{1245}-e^{1236}-e^{3678}-e^{4578}+e^{1568}+e^{2567}
\cr
e^{2348}-e^{1347}-e^{1246}+e^{1235}-e^{4678}+e^{3578}+e^{2568}-e^{1567}}
\label{vanishcayl}
\eeq

The Cayley calibration is a degree four calibration in an
eight-manifold $M$. Let $X$ be a generalized Cayley submanifold
whose tangent directions are spanned by $\{e^a; a=1,\dots,4\}$ and
normal directions by $\{e^i; i=5,\dots, 8\}$. Then, by the same
reasoning used for the generalized associative deformations, the
condition $\cL_V \tau |_X =0$ implies that
\beq
{\cal D}V_i:=\sum_{a,j}{t}^a_{ij}
\tilde\nabla_{a} V^j=0\ ,
\label{stuffia}
\eeq
where $\{t^a; a=1,\dots 4\}=\{ 1,\Omega_r;
r=1,2,3\}$, and $\{\Omega_r; r=1,2,3\}$ is a basis of constant
anti-self-dual 2-forms in $\bR^4$ spanned by the directions
$5,6,7,8$; such a basis has been defined in section nine for the
$G_2$ calibrations.
 The operator ${\cal D}$ is elliptic.  So if the moduli space
 exists, it is finite dimensional.
It is expected that for generic generalized  Cayley cycles, the
dimension of the moduli space is the index of the operator ${\cal D}$.
The index of this operator ${\cal D}$ is the same as that computed for
the standard Cayley calibrations because the principal symbol is
the same. It has been found \cite{joyceb} that ${\rm
ind}({\cal D})=\sigma(X)-{1\over 2}\chi(X)-{1\over2} [X]\cdot [X]$,
 where $\sigma(X)$, $\chi(X)$ and $([X]\cdot [X])$ is the signature,
Euler number and self-intersection of the Cayley calibration $X$.

\section{Examples}

\subsection{A group manifold example}

Let  $M=S^3\times \tilde S^3\times S^1\times \tilde S^1$
equipped with the left invariant
metric
\beq
g=\sum_{i=1}^3 (\s^i)^2+\sum_{i=1}^3(\tilde\s^i)^2+(\s^0)^2
+(\tilde\s^0)^2\
\label{stuffib}
\eeq
and the left invariant self-dual four-form
\beq
\Phi=\s^{0123}+(\s^{01}-\s^{23})\wedge (\tilde\s^{01}
-\tilde\s^{23})+(\s^{02}
+\s^{13})\wedge (\tilde\s^{02}+\tilde\s^{13})+
(\s^{03}-\s^{12})\wedge (\tilde\s^{03}-\tilde\s^{12})+\tilde \s^{0123}
\label{stuffic}
\eeq
where $\{\s^i; i=1,2,3\}$ and $\{\tilde\s^i; i=1,2,3\}$ are the
left-invariant one-forms of $S^3$ and $\tilde S^3$, respectively,
and $\s^0$ and $\tilde \s^0$ is the bi-invariant one-form of $S^1$
and $\tilde S^1$, respectively. Both the metric and self-dual
four-form are parallel with respect to the connection associated
with the left-action. With an appropriate choice of orientation of
$M$ both submanifolds $S^1\times S^3$ and $\tilde S^1\times \tilde
S^3$ are generalized Cayley calibrations. We shall focus on the
investigation of the moduli space of the $S^1\times S^3$ calibration;
the study of the moduli space of $\tilde S^1\times \tilde S^3$ is
similar. Observe that all the submanifolds $S^1\times S^3\times
\{p\}$, where $p\in \tilde S^1\times \tilde S^3$ are Cayley
calibrations. Therefore the dimension of the moduli space is at
least four. In fact the dimension of the moduli space is exactly
four. To see this one uses the fact that $\tilde \nabla=\partial$
acting on the normal vector fields of the calibration. Then the
result follows from a partial integration argument as in the group
manifold example in the $G_2$ case. The moduli space of the
$S^1\times S^3$ Cayley calibration is $\tilde S^1\times \tilde
S^3$. So $M$ is a family of generalized  calibrations for which
the fibers and the base are calibrated.

This example is a special case of a larger class of examples which
can be constructed by taking a seven-dimensional manifold
$(N,g,\psi)$ with a $G_2$ structure which admits a associative
submanifold $X$. Then the manifold $M=N\times S^1$ with
$ds^2(M)=ds^2(S^1)+ds^2(N)$ and $\Phi= e^0\wedge \psi+*\psi$  is a
$Spin(7)$ manifold; $*$ is the Hodge operation in $N$ and $e^0$ is
the invariant one-form along $S^1$. In addition $S^1\times X$ is a
generalized Cayley calibration.

\section{The $\partial\bar\partial$-lemma and some useful formulae}

The main tool that we shall use for the investigation of hermitian
manifolds with trivial canonical bundle
 is the $\partial\bar\partial$-lemma.
This can be stated as follows. Let
$(M,g,J)$ be a hermitian manifold, and $\phi$ and $\tilde \phi$ be
 two closed (p,q)-forms. Locally one can always
write
\begin{equation}
\label{ddbar} \phi=\tilde \phi+\partial \bar \partial \psi\ ,
\end{equation}
for some locally defined  (p-1, q-1)-form $\psi$. The
$\partial\bar\partial$-lemma states that if $\phi$ and $\tilde
\phi$ represent the same class in the Dolbeault cohomology, then
(\ref{ddbar}) is valid for some (p-1, q-1)-form on $M$.

Let $(M,g,J)$ be a $2n$-dimensional ($n>1$) Hermitian manifold
with complex structure $J$ and compatible Riemannian metric $g$.
Denote the K{\"a}hler form by $\Omega$.  The definitions $JX$ and
$J \alpha$, for $X$ a vector field and $\alpha$ a one-form, are
$(JX)^i=J^i{}_j X^j, \quad (J\alpha)^i=- (\alpha \circ
J)^i=-\alpha_k J^{ki} $ respectively. The Lee form $\theta$ is
defined by
\beq
 \theta=\dd\Omega \circ J\qquad \theta_i=-(\nabla^g)^k
\Omega_{kj} J^j{}_i
\label{stuffid}
\eeq
 where $\dd$ is the adjoint of $d$ and
$\nabla^g$ is the Levi-Civita connection of the metric $g$.
 Equivalently $\dd\Omega=J\theta. $ If the
 Lee form $\theta=0$ then the hermitian
manifold is said to be balanced. Balanced hermitian manifolds are
studied in \cite{Mi,Gau,AB1,AB2,GI1,GI2}.

The Bismut connection $\nabla^b$ and the Chern connection
$\nabla^c$ are given by
\begin{equation}\label{1}
g(\nabla^b_X Y,Z) = g(\nabla^g_X Y,Z) + \frac{1}{2} d^c \Omega
(X,Y,Z)\ ,
\end{equation}
\begin{equation}\label{1dop}
g(\nabla^c_X Y,Z) = g(\nabla^g_X Y,Z) + \frac{1}{2} d \Omega
(JX,Y,Z)\ ,
\end{equation}
respectively. Recall that $d^c = i(\overline {\partial} -
\partial )$. In particular, $d^c \Omega (X,Y,Z) = - d \Omega
(JX,JY,JZ)$.

Let $\rho^{\it b}$ and $\rho^{\it c}$ be the Ricci forms of the
Bismut and Chern connections respectively. Then it was shown in
\cite{AI1} that
\begin{equation}\label{ricc1}
\rho^{\it c} =\rho^{\it b} + d(J\theta)\ .
\end{equation}

In complex coordinates $\{z^{\alpha}\},
\alpha = 1,...,n$, we have the following formulae:

Let $\cG$ be the Chern connection. Then
\beq
(\cG)_{\alpha\beta}^{\delta}=g^{\bar \gamma\delta}
\partial_\alpha g_{\beta \bar \gamma}
\label{stuffie}
\eeq
and so
\beq
(\cG)_\alpha=(\cG)_{\alpha\beta}^{\beta}=
\partial_\alpha(\log(\det(g))\ .
\label{stuffih}
\eeq
The Lee form is then
\beq\label{ad2} \theta_\alpha =
(\cG)_{\alpha\beta}^{\beta}-(\cG)_{\beta\alpha}^{\beta}= g^{\beta
\bar\gamma}\left(\partial_\alpha g_{\beta \bar
\gamma}-\partial_\beta g_{\alpha \bar
\gamma}\right)=\partial_{\alpha}(\log(\det(g))-g^{\beta
\bar\gamma}\partial_\beta g_{\alpha \bar \gamma}\ .
\eeq
In terms of the Chern connection, the Lee form is
\beq\label{ad3}
(d\theta)_{\alpha\gamma}=\partial_{\alpha}(\cG)_{\beta\gamma}^{\beta}-
\partial_{\gamma}(\cG)_{\beta\alpha}^{\beta} ~.
\eeq
The Ricci form of the Chern connection is $\rho^{\it c}$
\beq\label{ad4}
 i\rho^{\it c}_{\bar{\beta}\alpha} = \partial_{\bar
\beta}D_{\alpha\sigma}^{\sigma}=
\partial_{\bar
\beta}\partial_{\alpha}(\log(\det(g))=\partial_{\bar\beta}
\left(g^{\bar \gamma\delta} \partial_{\alpha}g_{\delta \bar
\gamma}\right) \ .
\eeq

The (1,1)-part of formula (\ref{ricc1}) can be
written in the following way \cite{IP1}
$$
i\rho^{\it b}_{\bar{\beta}\alpha} =i\rho^{\it c}_{\bar{\beta}
\alpha}-\left(
\partial_{\bar \beta}\theta_{\alpha}+ \partial_{
\alpha}\theta_{\bar\beta}\right)~.
$$
Using (\ref{ad2}), (\ref{ad4}) we obtain that (\ref{ricc1})
 is equivalent
to the following two formulae
 \beq\label{ad5}
 i\rho^{\it b}_{\bar{\beta}\alpha}= \partial_{\bar{\beta}}\left(
 g^{\sigma
\bar\gamma}\partial_\sigma g_{\alpha \bar \gamma}\right) -
\partial_{\alpha}\left(
 g^{\sigma
\bar\gamma}\partial_{\bar{\beta}} g_{\sigma \bar \gamma}\right)+
\partial_{\alpha}\left(
 g^{\sigma
\bar\gamma}\partial_{\bar{\gamma}} g_{\sigma \bar \beta}\right)\ ,
\eeq \beq\label{ad6}
i\rho^{\it b}_{\beta\alpha}=(d\theta)_{\beta\alpha}=
\partial_{\beta}\theta_{\alpha}
- \partial_{\alpha}\theta_{\beta}\ .
\eeq

\section{Chern Connections with holonomy contained in $SU(n)$}

Hermitian manifolds which admit a Chern connection with holonomy
contained in $SU(n)$ necessarily admit a holomorphic (n,0)-form.
Since the existence of holomorphic (n,0)-forms depends only on the
choice of complex structure, this can be used to show whether a
complex manifold admits a hermitian structure for which the
associated Chern connection has holonomy contained in $SU(n)$. To
find whether a certain complex manifold admits (n,0)-holomorphic
forms one can use the Kodaira-type vanishing theorems
\cite{KW,Gau} together with the results in \cite{IP1}.

For complex manifolds satisfying the $\partial\bar{\partial}$-lemma,
we have the following:

\begin{theorem}\label{th2}

Let $(M,J)$ be a 2n-dimensional compact  connected complex
non-K{\"a}hler manifold with vanishing first Chern class,
$c_1(M,J)=0$. Suppose $(M,J)$ satisfies the
$\partial\bar{\partial}$-lemma. Then \item i) There exists a
Hermitian structure such that the restricted holonomy of the Chern
connection is contained in SU(n); \item ii) The Hodge number
$h^{n,0}\leq 1$. If $M$ is simply connected then $h^{n,0}=1$.

\end{theorem}

{\it Proof:} The Ricci form $\rho^{\it c}$ of the Chern  connection of
any hermitian structure $(g,J)$ represents the first Chern class
of the manifold. Therefore since $c_1(M,J)=0$, $\rho^{\it c}$ is
exact. Because $\rho^{\it c}$ is also a (1,1)-form, applying
 the $\partial\bar{\partial}$-lemma,  we find that
\beq
\rho^{\it c}=i\partial\bar{\partial}h~,
\label{stuffii}
\eeq
for some real function $h$ on
$M$.
Next we consider the manifold $M$ with hermitian
 structure $(M,\tilde g=e^{h/n}g,J)$.
Using (\ref{ad4}), we find that the Chern Ricci form
$\tilde \rho^{\it c}$ of the
new hermitian structure vanishes  because
\beq
i\tilde{\rho}_{\bar{\beta}\alpha}^{\it c}=\partial_{\bar{\beta}}
\partial_{\alpha}h + i\rho_{\bar{\beta}\alpha}^{\it c} =0.
\label{stuffij}
\eeq
Hence, the restricted holonomy of the Chern connection of $(M,
\tilde g J)$ is contained in $SU(n)$. This proves (i).

To show (ii),  the  Gauduchon plurigenera theorem \cite{Gau}
implies $h^{n,0}(M,J)\leq 1$ since the function $\tilde u
=trace(\tilde{\rho}^{\it c})=0$. If $M$ is simply connected, since
the holonomy of the Chern connection of $(M, \tilde g J)$ is
contained in $SU(n)$, there is a parallel (n,0)-form. A parallel
(n,0)-form with respect to the Chern connection is necessarily
holomorphic. Hence, $h^{n,0}=1$. This proves  (ii).
\hfill {\bf Q.E.D.}

\begin{corollary}

On $k\geq 2$-copies of $S^3\times S^3$ there exists a hermitian
structure such that the holonomy of the Chern connection is
contained in SU(3). In the conformal class of  any hermitian
structure there exists a unique (up to homothety) one with
$Hol(\nabla^c)\subseteq SU(3)$.

\end{corollary}

\section{Bismut Connections with holonomy contained in $SU(n)$}

Consider the following lemma

\begin{lemma}\label{lem1}

The Bismut Ricci forms $\rho^{\it b},\tilde{\rho}^{\it b}$
of two conformally
equivalent hermitian structures $(M,g,J)$ and $(M, \tilde g=e^fg,J)$ are
related by
\beq
i\tilde{\rho}_{\bar{\beta}\alpha}=i\rho_{\bar{\beta}\alpha}+
(2-n)
\partial_{\bar{\beta}}\partial_{\alpha}f \ ;
\quad \tilde{\rho}_{\beta\alpha}
=\rho_{\beta\alpha} \ .
\label{stuffkb}
\eeq

\end{lemma}

{\it Proof:}. It follows by straightforward
calculations from (\ref{ad5}) and (\ref{ad6}).
\hfill {\bf
Q.E.D.}

\begin{theorem}\label{th3}

Let $(M,J)$ be a 2n-dimensional compact complex non-K{\"a}hler
manifold with vanishing first Chern class, $c_1(M,J)=0$. Suppose
$(M,J)$ satisfies the $\partial\bar{\partial}$-lemma and
 that there exists a hermitian structure
$(g,J)$ such that $d\theta$ is a (1,1)-form. Then there exists
another unique (up to homothety) conformal hermitian structure
$(M,\bar g=e^fg,J)$ such that the restricted holonomy of the
associated Bismut connection is contained in $SU(n)$ provided
$n\geq 3$.

\end{theorem}

{\it Proof:}
Let  $(M,g,J)$ be the hermitian structure with
 $d\theta$ a (1,1)-form. Using (\ref{ad6}), we find
 that   the  Ricci form $\rho^{\it b}$
of the Bismut connection is a (1,1)-form. Therefore it is an
exact (1,1)-form since $\rho^{\it b}$ represents the first Chern
class which is zero. Applying the $\partial\bar{\partial}$-lemma
we can write $\rho^{\it b}=i\partial\bar{\partial}f$, for some
real function $f$ on $M$. Next using lemma~\ref{lem1}, it is
straightforward to observe that the Ricci form $\bar{\rho}^{\it
b}$ of the Bismut connection of the hermitian structure $(M,\bar
g=e^{f/2-n}g,J)$ vanishes.  Thus, the holonomy of the Bismut
connection of $(M,\bar g,J)$ is contained in SU(n). The uniqueness
follows since on a compact hermitian manifold the equation
$g^{\bar{\beta}\alpha}\partial_{\bar{\beta}}\partial_{\alpha}f=0$
has only constant solutions. This completes the proof. \hfill {\bf
Q.E.D.}

We remark that if the Lee form $\theta$ is an exact form i.e. the
structure is conformally balanced, then the above theorem applies.
So we have the following corollary:

\begin{corollary}\label{co2}

Let $(M,g,J)$ be a 2n-dimensional compact complex balanced
non-K{\"a}hler manifold with vanishing first Chern class,
$c_1(M,J)=0$. Suppose $(M,J)$ satisfies the
$\partial\bar{\partial}$-lemma. Then there exists another
conformal hermitian structure $(M,\bar g=e^fg,J)$, unique up to
homothety, such that the restricted holonomy of the associated
Bismut connection is contained in $SU(n)$ provided $n\geq 3$.

\end{corollary}

A Moishezon manifold is a compact complex manifold which is
bimeromorphic to a projective variety. Any Moishezon
manifold satisfies the $\partial\bar\partial$-lemma by a
 result of Deligne. Alessandrini and Bassanelli proved
in  \cite{AB2}, Corollary 4.6 that every Moishezon manifold is
balanced. Therefore from the above Corollary~\ref{co2} we have the
following:

\begin{corollary}\label{co3}

Every Moishezon manifold of complex dimension $n$, $n\geq 3$,
 with vanishing first Chern class admits a hermitian
structure for which both the Chern and Bismut connections
have restricted holonomy contained in $SU(n)$.

\end{corollary}

In particular Corollary~\ref{co3} implies that the Moishezon
manifolds with vanishing first Chern admit a holomorphic
(n,0)-form. Using the result shown in \cite{IP1},
 one concludes that the torsion of the Bismut
connection of Moishezon manifolds which are not K{\"a}hler and
 ${\rm Hol}(\nabla^b)\subseteq SU(n)$ is {\it not} closed.
Therefore such manifolds may have applications in heterotic strings.

\section{Proof of Theorem~\ref{t1}}

Let $M$ denote a connected sum of $k\geq 2$-copies of $S^3\times
S^3$. $M$ is cohomologically K{\"a}hler,  $b_1(M)=h^{1,0}=h^{0,1}=0$
and $h^{3,0}=h^{0,3}=1$ \cite{LuT}. According to Theorem 4.10 in
\cite{IP1}, if there exists a hermitian structure such that the
restricted holonomy $Hol(\nabla^b) \subset SU(3)$, then the
structure is conformally balanced i.e. its Gauduchon metric is
balanced. Conversely, the existence of a balanced Hermitian
structure on $M$ leads to the existence of hermitian structure in
the same conformal class with $Hol(\nabla^b) \subset SU(3)$ by
Theorem~\ref{th3}.

The existence of a balanced hermitian structure on a
2n-dimensional compact complex manifold $(M,J)$ has an intrinsic
characterization, namely it can be expressed in terms of positive
currents by a theorem of Michelsohn \cite{Mi} which states that a
smooth compact complex 2n-dimensional manifold admits a balanced
structure if and only if it is homologically balanced.

We recall that the space of real currents of degree (n-k,n-k) is
the dual space of $\Lambda^{k,k}(M)_R$, i.e. real (n-k,n-k)-forms
with distribution coefficients. A compact complex manifold is
balanced if and only if there is no positive current $T$ of degree
(1,1) which is the component of a boundary (i.e. if
$T=\bar{\partial}S+\partial\bar{S}$ and $T>0$ then $T=0$
\cite{Mi}). This result has an expression in terms of Aeppli group
$V^{1,1}(M)_R$. The real (1,1)-Aeppli group is defined as
 \beq
 V^{1,1}(M)_R=\frac{Ker(i\partial\bar{\partial}:
\Lambda^{1,1}(M)_R\to
 \Lambda^{2,2}(M)_R)}{(\partial\Lambda^{0,1}(M)+\bar{\partial}
\Lambda^{1,0}(M))_R}\
 \label{stuffkc}
 \eeq

The Michelsohn theorem can be rewritten \cite{AB2}: $M$ is
balanced if and only if every non-zero positive
$\partial\bar{\partial}$-closed current of degree (1,1) represents
a non-zero class in $V^{1,1}(M)_R$.

Now, it is clear that any compact cohomologically K{\"a}hler complex
manifold is balanced. Then, the above mentioned  result of Deligne
completes the proof of Theorem~\ref{t1} . \hfill {\bf Q.E.D.}

\begin{corollary}\label{co4}

Any compact 2n-dimensional $(n>2)$ complex manifold with vanishing
first Chern class which is cohomologically K{\"a}hler admits a
hermitian structure with vanishing Ricci form of the Bismut
connection.

\end{corollary}

\section{Examples}

Here we shall give examples of Bismut connections with holonomy
$SU(n)$ which can be thought of as generalizations of Calabi-Yau
manifolds.

{\bf Example 1}
Consider the $U(n)$ invariant metric
\beq
ds^2=(A(r^2) \delta_{\alpha \bar \beta}+ B(r^2)\bar z_{\alpha}
z_{\bar \beta}) dz^\alpha d\bar
z^{\bar \beta}
\eeq
where $ \bar z_{\alpha}=\delta_{\alpha\bar\beta}
 \bar z^{\bar \beta}$ and
$ z_{\bar\alpha}=\delta_{\bar\alpha\beta}  z^{\beta}$ and
$r^2=\delta_{\alpha\bar\beta} z^\alpha
z^{\bar \beta}$.
In this case it can be easily seen that the connection
 of the canonical
bundle is $\omega_\alpha= i\bar z_\alpha f(r^2)$,
where
\beq
f= (n-1) A^{-1} (2B-A')+ \big({\rm log}(A+r^2 B)\big)'
\eeq
where prime denotes differentiation with respect to $r^2$.
 The condition that
$d\omega=0$, necessary for the holonomy to be contained
 in $SU(n)$, implies that
\beq
f=0 \ .
\eeq
We remark that the K{\"a}hler case corresponds to taking
 $B=A'$. In this case
the solutions produce the Calabi-Yau metrics due to Calabi.

{\bf Example 2}. A compact example of a hermitian manifold for
which its Bismut connection has holonomy  SU(3) is as follows.
Consider the
 complex Heisenberg group
\beq
G=
\left\{ \left(
\begin{array}{lll}
1 & z_1  & z_3 \\
0 & 1    & z_2 \\
0 & 0    & 1   \\
\end{array}
\right) \quad z_1,z_2,z_3 \in {\bf C} \right\},
\label{examplestuffa}
\eeq
 with multiplication. The complex Iwasawa manifold is the compact
quotient space $ M=G/\Gamma $ formed from the right cosets of the
discrete group $\Gamma $ given by the matrices whose entries
$z_1,z_2,z_3$ are Gaussian integers. The $1$-forms  $dz_1,\quad
dz_2,\quad dz_3-z_1dz_2$ are left invariant by $G$ and  by
$\Gamma$. These 1-forms pass to the quotient $M$.  We denote by
$\alpha_1,\alpha_2,\alpha_3$ the corresponding 1-forms on $M$,
respectively. Consider the Hermitian manifold $(M,g,J)$, where $J$
is the natural complex structure on $M$ arising from the complex
coordinates $z_1,z_2,z_3$ on $G$ and the metric $g$ is determined
by $g=\sum_{i=1}^3\alpha_i\otimes \bar{\alpha }_i$. The Chern
connection $D$ is determined by the conditions that the 1-forms
$\alpha_1,\alpha_2,\alpha_3$ are parallel. The torsion tensor of
$D$ is given by
$C(\alpha_i^{\#},\alpha_j^{\#})=-[\alpha_i^{\#},\alpha_j^{\#}],
\quad i,j=1,2,3, $ where $\alpha_i^{\#}$ is the vector field
corresponding to $\alpha_i$ via $g$. The only nonzero term is
$C(\alpha_1^{\#},\alpha_2^{\#})=-\alpha_3^{\#}$ and its complex
conjugate. Thus, the space $(M,g,J)$ is a compact balanced
Hermitian (non K{\"a}hler) manifold with a flat Chern connection
and automatically the holonomy group of its Bismut connection is
contained in SU(3) by formula (\ref{ricc1}). The (0,3)-form
$\psi = \alpha_1\wedge\alpha_2\wedge\alpha_3$ is parallel with
respect to both Chern and Bismut connections. Let
$e_1,e_2,e_3,Je_1,Je_2,Je_3$ be a real basis determined by
$\alpha_j^{\#}=e_j-\sqrt{-1}Je_j$. Then the real Iwasawa manifold
$X$ determined by $e_1,e_2,e_3$ is a SAS calibration with respect
to $Re \psi$. It admits moduli since $V=-Je_3$ is a SAS
deformation. Indeed, $U=JV=e_3$ is a Killing vector field on $M$
and therefore it is also  holomorphic by results in \cite{GI2}
since the Chern connection is flat. Thus, $e_3$ is a SAS
deformation of $X$ in $M$.

\end{document}